\documentclass[preprint,12pt]{elsarticle}

\usepackage{amssymb}
\usepackage{amsmath}
\usepackage{graphicx}%
\usepackage{multirow}%
\usepackage{amsmath,amssymb,amsfonts}%
\usepackage{amsthm}%
\usepackage{mathrsfs}%
\usepackage[title]{appendix}%
\usepackage{xcolor}%
\usepackage{textcomp}%
\usepackage{manyfoot}%
\usepackage{booktabs}%
\usepackage{algorithm}%
\usepackage{algorithmicx}%
\usepackage{algpseudocode}%
\usepackage{listings}%
\usepackage{eurosym}%
\usepackage{tabularx}%
\usepackage{empheq}%
\usepackage{subcaption}%
\usepackage{xurl}

\newcommand{\soc}{\text{soc}}
\newcommand{\SOC}{\text{SoC}}

\newcommand{\R}{\mathbb{R}}

\journal{EURO Journal on Computational Optimization}

\begin{document}

\begin{frontmatter}



\title{Optimal Participation of Energy Communities in Electricity Markets under Uncertainty. A Multi-Stage Stochastic Programming Approach}


\author[UPC]{Albert Sol\`a Vilalta} 
\ead{albert.sola.vilalta@upc.edu}

\author[UPC]{Ignasi Ma\~{n}\'e Bosch}

\author[UPC]{F.- Javier Heredia}
\ead{f.javier.heredia@upc.edu}

\affiliation[UPC]{organization={Departament d'Estadística i Investigació Operativa, Universitat Politècnica de Catalunya},
            addressline={C/ Jordi Girona, 31}, 
            city={Barcelona},
            postcode={08034},
            country={Spain}}

\begin{abstract}
We propose a multi-stage stochastic programming model for the optimal participation of energy communities in electricity markets. The multi-stage aspect captures the different times at which variable renewable generation and electricity prices are observed. This results in large-scale optimization problem instances containing large scenario trees with 34 stages, to which scenario reduction techniques are applied. Case studies with real data are discussed to analyse proposed regulatory frameworks in Europe. The added value of considering stochasticity is also analysed.
\end{abstract}




\end{frontmatter}



\section{Introduction} \label{SEC:Introduction}

In 2018 and 2019, the European Union passed legislation to create energy communities. They are new legal entities that empower citizens, small businesses, and local authorities to produce, manage, and consume their own energy \cite{EU2018/2001, EU2019/944}. Their implementation aims to increase energy efficiency, increase renewable generation, contribute towards the electrification of the heating and transport sectors, and encourage citizens' behavioural change. These are four of the pillars the International Energy Agency has identified on their global pathway to net zero CO$_2$ emissions \cite{IEA21}.

In this work, we focus on the electricity sector and, in particular, try to address the following research question: How should an energy community optimally satisfy its electricity needs and make the most of its flexibility in wholesale electricity markets? The proposed approach can be adapted to model any entity with both flexible electricity demand and variable renewable generation that is able to participate in electricity markets, as is the case for aggregators.


The electricity markets we consider in this work are the day-ahead market, secondary reserve market, and intraday markets, which are complemented with an imbalance settlement. Together, they represent the core structure of electricity markets in those regions where electricity supply is liberalised \cite{Wolak21, WorldBank22}. 

Within this framework, the decision-making problem faced by the energy community is affected by multiple uncertainties, which can be classified into two different groups. First, it is affected by uncertainty in their own variable renewable generation, typically solar and/or wind. Second, it is affected by uncertainty in electricity prices in all markets it wishes to participate because the decision on which bid to submit is made before the markets are cleared. To cope with all these uncertainties, we propose a multi-stage stochastic programming model  \cite{Pflug14} that maximizes the expected social welfare of an energy community from satisfying its internal electricity demand and participating in electricity markets.



In a stochastic programming approach, to make the resulting optimization problem tractable, all random variables are assumed to be discrete. To capture the multi-stage structure of the problem, a scenario tree is built following a two-phase method adapted from \cite{Cuadrado20} to generate the required scenario trees from historical data, where the branching of the tree represents moving from one stage to the next one. 

The main contribution of this work is the introduction of a novel multi-stage stochastic programming model for the optimal participation of energy communities in the day-ahead, reserve and intraday electricity markets. The decision-making process is modelled in great detail, including a stage for every electricity market where the energy community might participate, and for every hour of the day, to accurately model the hourly behaviour of the energy community. It is based on previous work on the optimal participation of a virtual power plant, formed by a wind farm and a battery, in electricity markets \cite{Heredia18, Cuadrado20}.

In this regard, the novelties introduced in the model are the following:
\begin{enumerate}
    \item Introduction of combined buying and selling bids. These are bids that buy energy if the price is sufficiently low, and sell energy if the price is sufficiently high.
    
    \item Co-optimization of buying and selling bids. This allows the energy community to better cope with electricity price and variable renewable uncertainty.
    
    \item Novel modelling of aggregated flexible demand. It includes enough detail to limit the flexibility in the whole planning horizon and on specific periods, and is suitable to fit a complex stochastic programming environment.
    
    \item Generation of scenario trees for the stochastic programming model, extending the methodology in \cite{Cuadrado20, Heredia18} to include solar generation.
\end{enumerate}

After this introductory section, the remainder of this article is organized as follows. Section \ref{SEC:LiteratureReview} reviews the literature on energy communities and flexible demand modelling, European electricity market regulation, and stochastic programming approaches for multi-market participation. Section \ref{SEC:TheDecision-MakingProcessUnderUncertainty} describes the decision-making process faced by an energy community that wishes to participate in electricity markets. Section \ref{SEC:GenerationOfScenariosTrees} describes the methodology to generate scenario trees and Section \ref{SEC:TheMulti-StageStochasticProgrammingModel} presents the multi-stage stochastic programming model. Section \ref{SEC:ModelOutputs} describes the main model outputs and Section \ref{SEC:AnalysisOfResults} presents a $1$-month case study and analyses the results. Finally, conclusions are drawn in Section \ref{SEC:Conclusions}.

\section{Literature Review} \label{SEC:LiteratureReview}

\subsection{Energy Communities}
\label{SUBSEC:EnergyCommunities}

Energy communities were introduced in European legislation through European Directives (EU) 2018/2001 \cite{EU2018/2001} and (EU) 2019/944 \cite{EU2019/944}, which have been recently amended by Directives (EU) 2023/2413 \cite{EU2023/2413} and (EU) 2024/1711 \cite{EU2024/1711}. These directives are transposed into national legislation, adapting them to the specific characteristics of every country. For example, Directive (EU) 2018/2001 was transposed into Spanish legislation through Real Decreto-ley 23/2020 \cite{BOE23/2020} and into Italian legislation through Decreto Legislativo n. 199 \cite{GU21}.

Energy communities have been often modelled as entities that have an internal energy demand, energy storage systems, and generation capabilities, typically in the form of variable renewable generation \cite{Manso-Burgos22, Brusco23, Vecchi24, Crowley25, Nour25}. Previous works in the literature on energy communities can be classified in three main streams: distribution of benefits and responsibilities between energy community members, participation in electricity markets and regulation analysis.

The distribution of benefits and responsibilities between energy community members consists on analysing how a service (e.g., flexibility) or commodity (e.g., electricity) provided by the energy community should be obtained from its members, and what remuneration should members get for their contribution. In this stream of literature, the members' interest to be part of the energy community is also often analysed. \cite{Manso-Burgos22} study several asset ownership options within an energy community and compare two sharing strategies: fixed and variable coefficients. \cite{Crowley25} propose a bilevel optimisation model to coordinate an energy community manager with its prosumer members for the provision of capacity limitation services. Finally, \cite{Garcia-Munoz25} analyse the ability of an energy community to provide flexibility services by aggregating its member's ability to provide such services.

Regarding the participation of energy communities in electricity markets, the focus is on the interactions between energy communities and the outside world. \cite{Brusco23} analyse the participation of energy communities in the day-ahead market and the provision of ancillary services to the power system using a deterministic model with simulated data. \cite{Orozco22} propose a multistage stochastic optimization model for day-ahead scheduling combined with an intraday rolling-horizon procedure. Other works in this stream include \cite{Long18} and \cite{Ferrucci25}, where the provision of grid ancillary services by an energy community is studied. 

In many aspects, this problem is similar to that faced by an aggregator of consumers and/or prosumers participating in electricity markets. \cite{Henriquez18} propose a two-stage stochastic bilevel optimization problem to determine the optimal participation of a demand response aggregator in day-ahead and real time markets. \cite{Iria18} present a two-stage stochastic programming model for the optimal participation of a prosumer aggregator in day-ahead and real-time markets, which is generalized in \cite{Iria19} to include tertiary reserve bids.

Finally, with the energy community regulation being in active development, several works have compared regulatory frameworks in different European countries and proposed improvements. \cite{Ines20, Haji-Bashi23} review and compare the regulatory frameworks for energy communities in different European member states. \cite{Sokolowski20} proposes and analyses several regulatory improvements. \cite{Ferrucci24} examine market structures, regulatory environments and technical challenges to identify barriers and potential solutions to integrate energy communities into electricity markets.

\subsection{Flexible Demand}

Active citizen participation through energy communities should increase flexibility in electricity demand. \cite{Barbero23} identify three different types of electricity demand, depending on their flexibility: fixed loads, shiftable loads, and curtailable loads. In this work, we focus on fixed and shiftable loads.

\cite{Henriquez18} use duration differentiated loads contracts and load curtailment contracts to model fixed loads, shiftable loads and curtailable loads. \cite{Iria18, Iria19} model fixed loads and shiftable loads, including great detail in some shiftable loads from thermostatical controls. 

While these works model flexible demand in greater detail than we do, they rely on modelling loads individually. This greatly increases the number of decision variables needed in the optimization problem, resulting in intractable optimization problems once the complete decision process we aim to tackle is included.

\subsection{European Electricity Market Regulation}

Electricity market regulation in Europe has been changing rapidly in recent years, with the final goal to better integrate electricity markets and facilitate cross-border trade. \cite{EU2019/943} set the European electricity market framework in which energy communities were first defined. We use this framework for our modelling, and in particular, the transposition to Spanish law described in \cite{CNMC19} and \cite{CNMC21}.

In recent months, we have seen significant regulatory changes to further enhance electricity market integration. In particular, Spain saw the introduction of the continuous intraday market, the reduction from $6$ to $3$ intraday market auctions and the move to $15$ minute periods in both intraday and balancing markets \cite{BOE24_IM, BOE24_balances}, which came into force in June 2024. Furthermore, the day-ahead market is expected to move to $15$ minute periods in late 2025 \cite{ENTSOE25}. We leave the adaptation of the model to the recent and coming regulatory changes for future work.

\subsection{Stochastic Programming for Multi-Market Participation}

The proposed multi-stage stochastic programming model for the optimal participation of energy communities in electricity markets follows the approach by \cite{Heredia18} and \cite{Cuadrado20} by extending the Virtual Power Plant consisting of a wind farm and a battery energy storage system (BESS) to an energy community. This required the inclusion of flexible demand and solar PV generation, the co-optimization of buying, selling and combined buying and selling bids, and the adaptation of the scenario trees generation procedure accordingly. 

Stochastic programming is widely used to optimize under uncertainty in many different problems, and energy communities are no exception. However, most of this works use simpler two-stage stochastic optimization models \cite{Henriquez18, Iria18, Iria19, Gu23, Tostado-Veliz23}. We have identified \cite{Orozco22} and \cite{Nemati25} as the works in the literature that are most similar to ours.

The model presented in \cite{Orozco22} is a multi-stage stochastic programming model for day-ahead and intraday market participation, but the hours are grouped in three groups of $8$ hours, resulting in $3$ stages. Therefore, while the stochastic programming model is more complex than most, it does not capture as much detail as our proposal on the decision-making process and does not model the reserve market.

\cite{Nemati25} propose a novel approach based on robust optimization to multi-market participation of a virtual power plant in day-ahead, reserve and intraday electricity markets. In this sense, the decision-making process is similar to \cite{Heredia18} and ours, but the approach is quite different. Moreover, their focus is on a virtual power plant and as a result there is no flexible demand.

\section{The Decision-Making Process under Uncertainty} \label{SEC:TheDecision-MakingProcessUnderUncertainty}


The decision-making problem faced by the energy community is affected by uncertainty in their own variable renewable generation and electricity market prices. Regarding variable renewable uncertainty, we assume that the energy community has solar and wind generation capacity. In the model, we aggregate them by type, as if there were a single solar farm and a single wind farm in the energy community. This is justified because all energy community assets should be relatively close to each other, and hence affected by the same weather conditions, which should lead to similar generation outputs. We assume that variable renewable generation at time $t$ is observed at the beginning of the time period $t$.

Regarding uncertainty on electricity market prices, we incorporate it in all the markets we consider. These are: the day-ahead market (DM), the secondary reserve market (RM) and the intraday markets (IM). We also incorporate uncertainty in the imbalance settlement (IB) prices. We follow closely the regulation of Spanish electricity markets before June 2024 \cite{CNMC19, CNMC21}, which is similar to other electricity market structures in Europe and other liberalized markets \cite{Wolak21, WorldBank22}. Next, we summarize the markets considered:

\begin{itemize}
    \item \textbf{Day-ahead market (DM)}: generates an hourly plan of the electricity exchanges for the next day. For electricity exchanges on day $D$, market participants submit bids at midday on day $D - 1$.

    \item \textbf{Secondary reserve market (RM)}: ancillary service whose objective is to maintain the generation-demand balance. Market participants submit availability capacity bids (MW). If the bid is accepted, they must be ready to provide that capacity within an activation time limit (5 minutes) for a certain period of time (15 minutes) at the instruction of the system operator.

    \item \textbf{Intraday markets (IM)}: allows market participants to adjust their day-ahead positions as delivery time approaches. \cite{CNMC21} establishes $6$ auction intraday electricity markets with different gate closures, and a continuous intraday electricity market. We consider only the $6$ auction intraday electricity markets.

    \item \textbf{Imbalance settlement (IB)}: compensates or penalizes market participants depending on whether their committed positions in day-ahead and intraday markets were achieved or not. It is done after delivery time, on day $D + 1$.
\end{itemize}

\begin{figure}[ht]
    \centering
    \includegraphics[width=\linewidth]{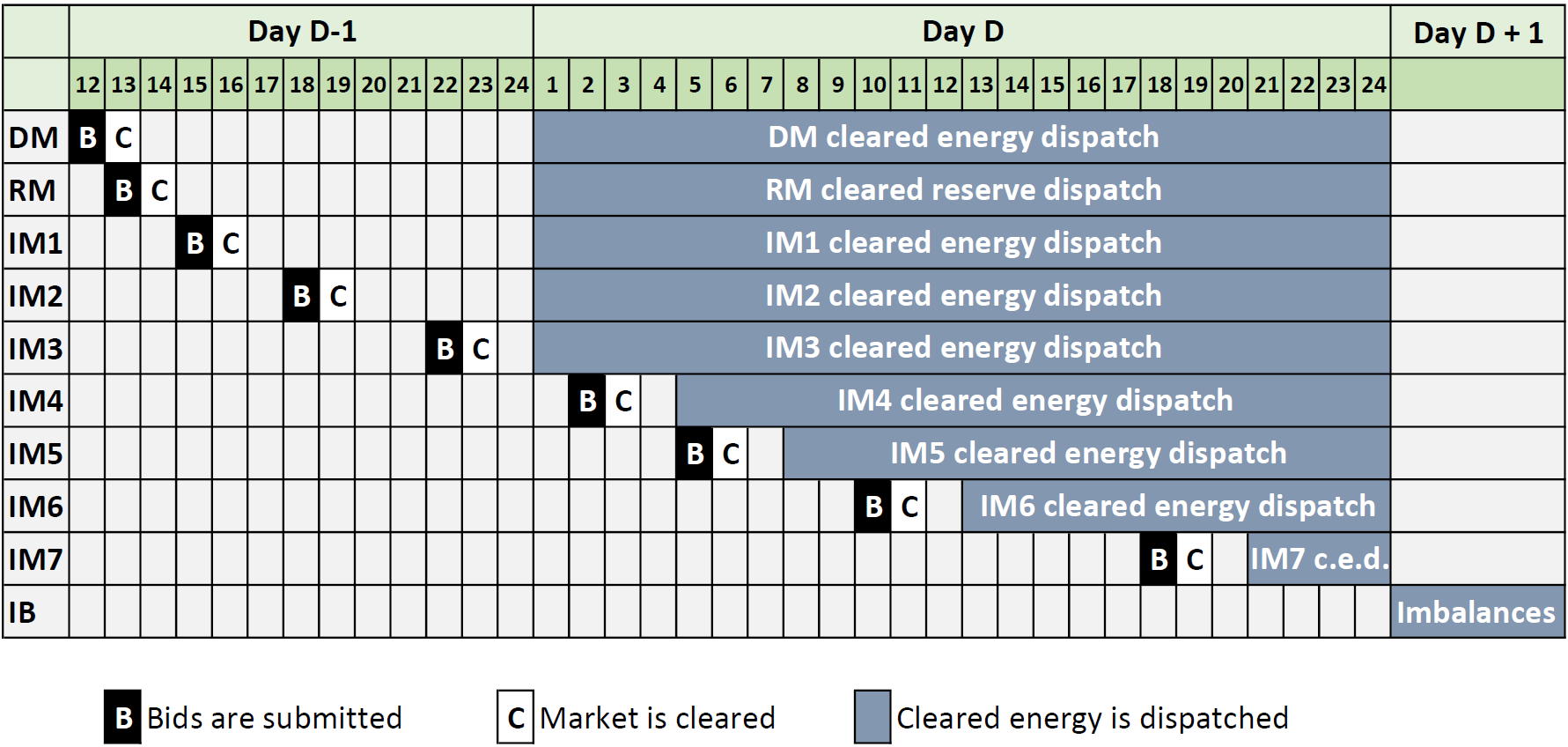}
    \caption{Timeline of the bid submission windows and market clearing processes of the markets considered.}
    \label{FIG:market_stages}
\end{figure}

Table \ref{TBL:RVconsidered} summarises the random variables considered and Figure \ref{FIG:market_stages} presents a timeline of the bid submission windows and market clearing processes for all markets considered. They correspond to the day-ahead, secondary reserve and intraday markets in Spain before June 2024. Observe that Figure \ref{FIG:market_stages} contains a seventh intraday market (IM7), which represents the first $4$ hours of the second intraday market (IM2) of day $D + 1$. This is because in the second intraday market of day $D + 1$, bids are accepted for the last $4$ hours of day $D + 1$ and the $24$ hours of day $D + 2$ \cite{CNMC21}. 

\begin{table}
\begin{tabular}{ l l }
	\hline
    $\lambda^D(\omega) \in \mathbb{R}^{24}$ & Day-ahead market clearing prices [\euro /MWh]. \\
	$\lambda^R(\omega) \in \mathbb{R}^{24}$ & Reserve market clearing prices [\euro /MW]. \\
    $\lambda_i^I(\omega) \in \mathbb{R}^{d_i}$ & Market clearing prices for intraday market $i \in \mathcal{I}$ [\euro /MWh], \\
    & $(d_1, \dots, d_7) = (24, 24, 20, 17, 13, 9, 4 )$. \\
    $W(\omega) \in \mathbb{R}^{24}$ & Wind generation [\euro /MWh]. \\
    $PV(\omega) \in \mathbb{R}^{24}$ & Solar PV generation [\euro /MWh]. \\
    $\lambda^{IB+}(\omega) \in \mathbb{R}^{24}$ & Positive imbalance prices [\euro /MWh]. \\
    $\lambda^{IB-}(\omega) \in \mathbb{R}^{24}$ & Negative imbalance prices [\euro /MWh]. \\
    \hline
\end{tabular}
\caption{Summary of the random variables considered.}
\label{TBL:RVconsidered}
\end{table}

\section{Generation of Scenarios Trees and Data} \label{SEC:GenerationOfScenariosTrees}

The stochastic process defined by all the uncertainties introduced in Section \ref{SEC:TheDecision-MakingProcessUnderUncertainty} is complex, and would result in an intractable optimization problem if incorporated in full detail. Hence, we discretise it and approximate it with a scenario tree, which preserves the timeline of the decision-making process. A scenario tree is a cycle-free, finite and directed graph with a single root for which the distance to all terminal nodes from the root is constant. In addition, the scenario tree carries probability valuations on nodes and arcs \cite{Pflug14}.

The stochastic process considered does not have a known probability distribution, but we can fit it with historical data. In liberalised electricity systems, market data is usually published by market operators, and Spain is no exception. OMIE, the Spanish and Portuguese market system operator, publishes hourly market data in its platform OMIEData \cite{OMIEData} for all markets we consider, and also the imbalance penalties and compensations. 

Regarding solar and wind generation data, it can be estimated using open-source tools like Renewables Ninja \cite{Pfenninger16, Staffell16, RenewablesNinja} or PVGIS \cite{PVGIS}. They both use meteorological data to estimate generation of a variable renewable installation in a given location. Renewables Ninja can provide solar and wind generation worldwide, whereas PVGIS can provide solar generation in Europe.

To generate scenario trees, we adapt the two-phase method from \cite{Heredia18, Cuadrado20} to the stochastic process described in Section \ref{SEC:TheDecision-MakingProcessUnderUncertainty}. Phase I consists on applying a Time Series Factor Analysis to describe the stochastic model with a reduced number of (unobserved) variables or factors \cite{Gilbert05, Gilbert06}. Vector autoregressive techniques are used to keep, to a significant extent, the correlation between the original random variables in the reduced model \cite{Sims80}. Finally, sampling from the obtained statistical models is performed using bootstrap techniques from \cite{Efron94,Munoz13}. At the end of Phase I, we obtain a scenario fan.

Phase II of the scenario tree generation procedure consists on applying a scenario reduction technique to reduce the size of the scenario fan while still representing the stochastic process in an accurate manner. To do so, we apply the Forward Tree Construction Algorithm developed in \cite{Dupacova03, Heitsch03}.


\section{The Multi-Stage Stochastic Programming Model} \label{SEC:TheMulti-StageStochasticProgrammingModel}


We consider an energy community with the following elements:
\begin{enumerate}
	\item Flexible Demand.
	\item Wind Farm.
	\item Solar Farm.
	\item Battery Energy Storage System (BESS).
\end{enumerate}

We present the complete mathematical formulation of the model in \ref{APP:MathematicalFormulation}.

\section{The Model's Optimal Solution} \label{SEC:ModelOutputs}

The optimal decisions included in the optimal solution of the model can be classified in three categories: optimal bid of the energy community to the day-ahead market, optimal price-accepting bid of the energy community to the reserve and intraday markets, and optimal behaviour of all elements of the energy community. They are described in Sections \ref{SUBSEC:OptimalBidToTheDayAheadMarket}, \ref{SUBSEC:OptimalPrice-AcceptingBidsToReserveAndIntradayMarkets} and \ref{SUBSEC:OptimalBehaviourOfAllElementsOfTheEnergyCommunity}.

\subsection{Optimal Bid to the Day-Ahead Market}
\label{SUBSEC:OptimalBidToTheDayAheadMarket}

The optimal bid of the energy community in the day-ahead market is a set of ordered quantity-price pairs $\{(q_j, p_j) \ | \ j \in \mathcal{J} \}$, where $\mathcal{J}$ is an index set over the quantity-price pairs or \emph{blocks} of the bid. Every pair $(q_j, p_j)$ indicates that the energy community is willing to buy (or sell) up to $q_i$ units of electricity [MWh] at the electricity market price $p_j$ [\euro /MWh]. A key feature of energy communities is that they both consume and generate electricity. Hence, they might be willing to sell or buy electricity, depending on the time period and/or electricity market price. To capture this behaviour, we use the convention that $q_j < 0$ represents a buying quantity and $q_j > 0$ a selling quantity.

In the literature, buying and selling bids are usually treated as separate curves, where both buying and selling quantities are nonnegative. Using this convention, the usual monotonicity bid requirement in electricity markets is that buying (resp. selling) bids must be monotonically decreasing (resp. increasing), i.e., for every $j, k \in \mathcal{J}$, if $q_j < q_k$, then $p_j > p_k$ (resp. $p_j < p_k$). Using our convention that $q_j < 0$ represents a buying quantity and $q_j > 0$ a selling quantity, the monotonicity bid requirement translates to bids being monotonically increasing. Observe that the monotonicity bid requirement aligns with the assumption that electricity market participants are economically rational, i.e., if electricity prices are higher, then their willingness to buy (resp. sell) is lower (resp. higher).

\begin{figure}[ht]
    \centering
    \includegraphics[width= \linewidth]{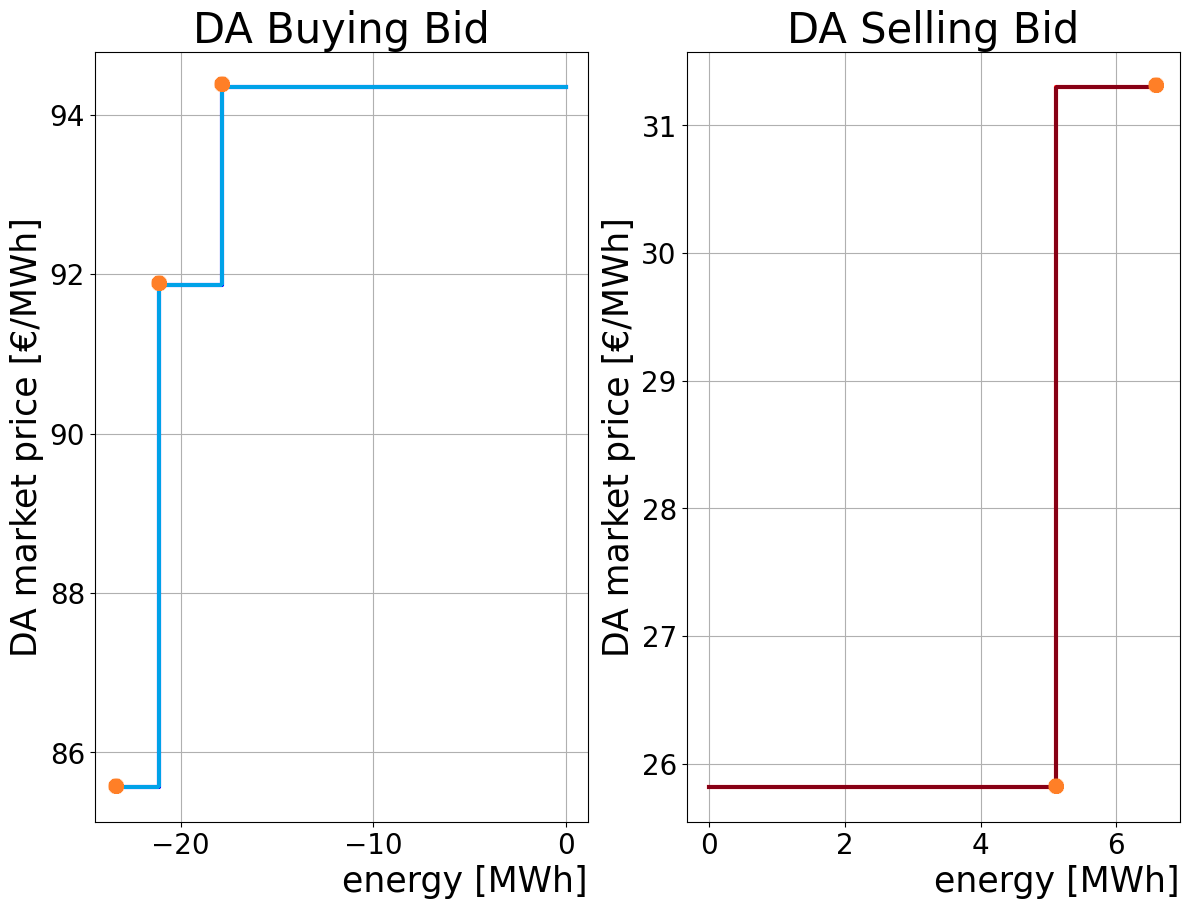}
    \caption{Day-ahead buying (left, blue) and selling (right, red) bid examples. The points in orange represent the quantity-price pairs $\{(q_j, p_j) \ | \ j \in \mathcal{J} \}$ that constitute the bid.}
    \label{FIG:DA_bid_examples}
\end{figure}

\begin{figure}[ht]
    \centering
    \includegraphics[width= \linewidth]{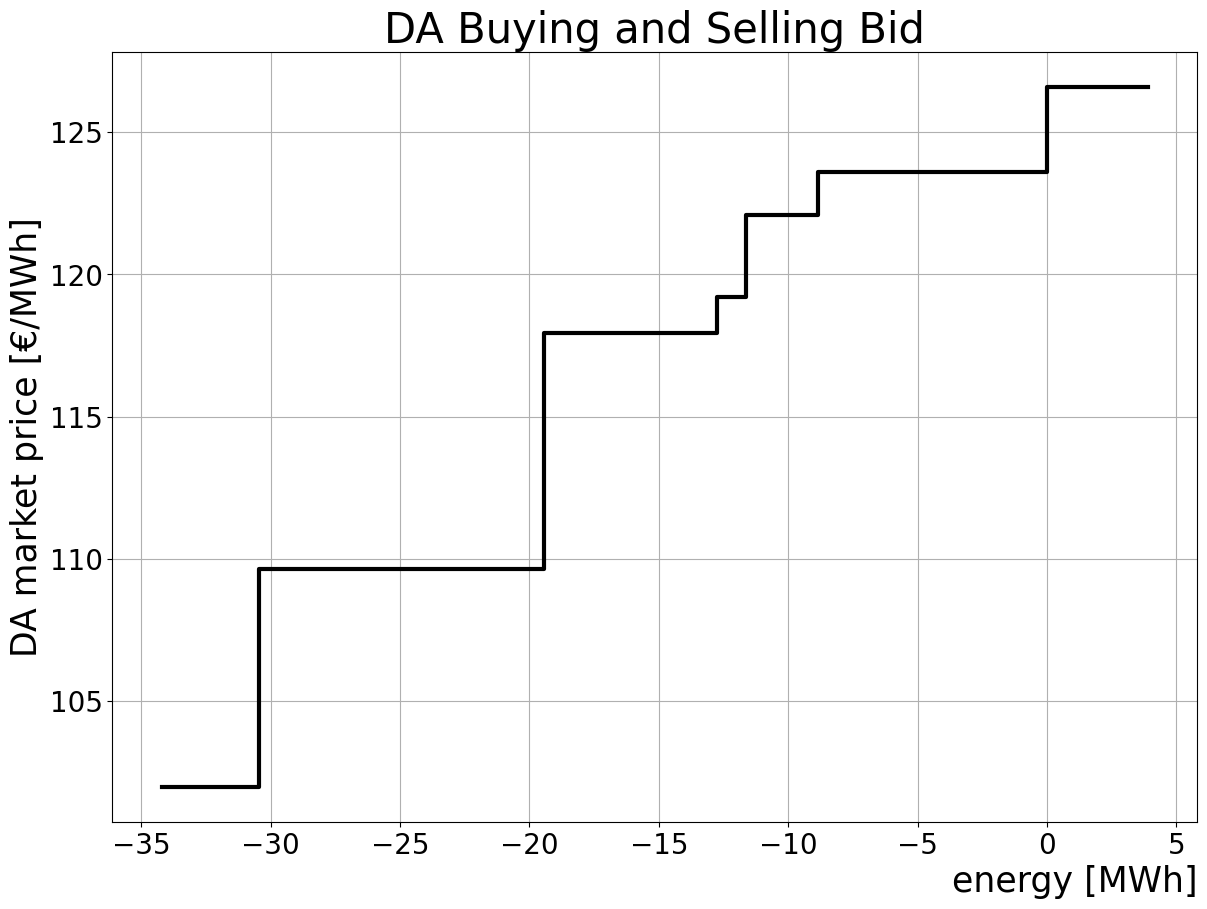}
    \caption{Combined buying and selling bidding curve.}
    \label{FIG:DA_bid_FTC_750_d009_t12}
\end{figure}

Figure \ref{FIG:DA_bid_examples} shows examples of buying and selling bids to the day-ahead market obtained by the model in the case study presented in Section \ref{SEC:AnalysisOfResults} at two different time periods. The quantity-price pairs $\{(q_j, p_j) \ | \ j \in \mathcal{J} \}$ are represented using orange points, which are joined to form the bidding curves. The bidding curves indicate the position of the energy community between quantity-price pairs. Figure \ref{FIG:DA_bid_FTC_750_d009_t12} shows an example of a combined buying and selling bidding obtained by the model, in which the energy community is willing to buy if electricity prices are below 124 \euro/MWh and sell if prices are above 126 \euro/MWh.

The bid's quantity-price pairs $(q_j, p_j)$ are obtained from the day-ahead energy matched variables $e^{DA+}_{t, \omega}$ and $e^{DA-}_{t, \omega}$, and the day-ahead market price parameters $\lambda^D_{t, \omega}$. We define the set of quantity-price pairs as $\mathcal{J} = \mathcal{C}_1$, the set of clusters at stage $1$. For every $j \in \mathcal{C}_1$, we denote by $\mathcal{C}_{1,j}$ the scenarios in the cluster $j$ and define the pair $(q_{j}, p_{j}) = (e^{DA+}_{t, \omega_j} + e^{DA-}_{t, \omega_j}, \lambda^D_{t, \omega_j})$ where $\omega_j$ is a scenario in $\mathcal{C}_{1, j}$. Note that the nonanticipativity constraints \eqref{EQ:NAC_e} guarantee that this definition is independent of the choice of $\omega_j \in \mathcal{C}_{1, j}$. Constraints \eqref{EQ:DA_sell_bid_LB} and \eqref{EQ:DA_buy_bid_LB} ensure that, at every scenario, only one of $e^{DA-}_{t, \omega_j}$ and $e^{DA+}_{t, \omega_j}$ is nonzero. Finally, the monotonicity constraints \eqref{EQ:DA_bid_mono_1} and \eqref{EQ:DA_bid_mono_2} guarantee that the resulting bid is monotonically increasing.

This approach has the advantage of reducing the amount of decision variables required to describe bidding curves, making the resulting optimization problem more tractable. The drawback is that the prices in the quantity-price pairs are fixed a priori by the scenarios, which reduces the flexibility of the model. In a stochastic programming environment where the electricity prices are determined by scenarios, however, this might be less of an issue. This is because the optimization model has no information other than the scenarios. We refer the reader to \cite{Corchero13} for a discussion on the optimality of this approach to obtain bidding curves in a stochastic programming environment, in the context of thermal generation.

\subsection{Optimal Price-Accepting Bids to Reserve and Intraday Markets}
\label{SUBSEC:OptimalPrice-AcceptingBidsToReserveAndIntradayMarkets}

We model reserve and intraday market participation differently than in day-ahead. Instead of aiming to build a bidding curve from the model as described in the previous paragraphs, we consider a single quantity for all scenarios at the respective clusters. This represents a single quantity for all reserve scenarios and a single quantity for every intraday market auction. They are often called \emph{price-accepting} bids because no matter the price at which the corresponding market clears, the bid quantity is the same.

Note that, since these decisions are taken at stages later than day-ahead, the corresponding variables ($r^U_{t, \omega}, r^D_{t, \omega}, e^{IM}_{i, t, \omega}$,...) can take different values at different scenarios, but this is because of the presence of different scenarios before the corresponding decision is made. This difference can be appreciated in Table \ref{TABLE:stages-rv-dv}, where reserve and intraday market variables ($r^U_{t, \omega}, r^D_{t, \omega}, e^{IM}_{i, t, \omega}$,...) are \emph{state variables} of the corresponding stages whereas day-ahead variables ($e^{DA+}_{t, \omega}, e^{DA-}_{t, \omega}, ie^{DA+}_{t, \omega}$) are \emph{recourse variables} of stage 1, the day-ahead stage. We follow the distinction between \emph{state variables} and \emph{recourse variables} in a multi-stage stochastic optimization problem from \cite{Escudero07, Cuadrado20}. In the multi-stage stochastic optimization model, this is introduced through the non anticipativity constraints, compare \eqref{EQ:NAC_r}, \eqref{EQ:NAC_rB}, \eqref{EQ:NAC_rFD} and \eqref{EQ:NAC_eIM} with \eqref{EQ:NAC_e} and \eqref{EQ:NAC_ie}.

The reason for not including the full detail of reserve and intraday market bids is to reduce the size of the problem and to keep it tractable. This is a reasonable approximation for two reasons. First, because we reduce the model's accuracy and flexibility at stages after the first stage, by which the first stage decisions should be less affected. Second, because, in a real-world situation, one could re-optimize with updated information once the corresponding decision has to be made. In this re-optimization process, it would make sense to include the full detail of bids as is currently done for day-ahead.

\subsection{Optimal Behaviour of All Elements of the Energy Community}
\label{SUBSEC:OptimalBehaviourOfAllElementsOfTheEnergyCommunity}

\begin{figure}[ht]
    \centering
    \includegraphics[width= \linewidth]{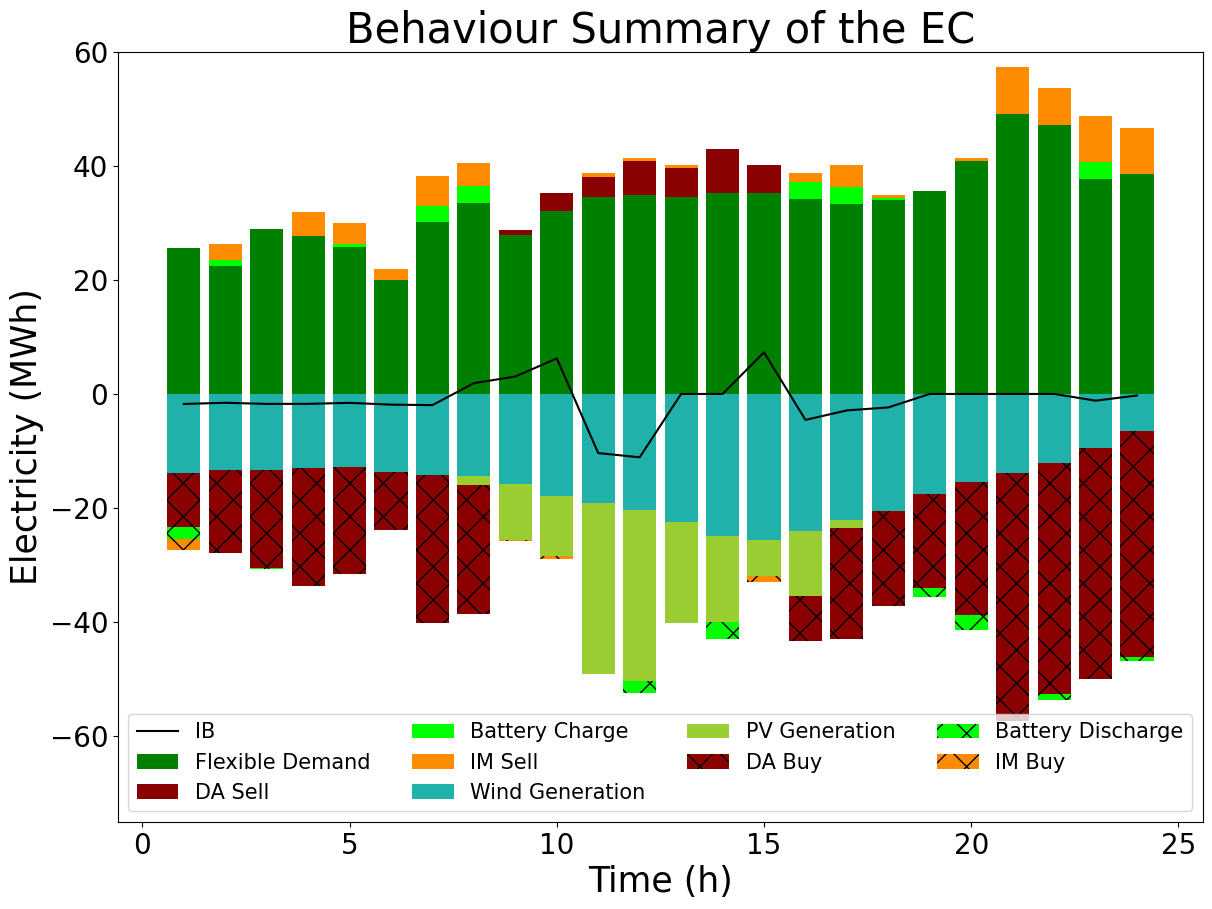}
    \caption{Behaviour summary of the energy community in a scenario during a whole day. Day-ahead position (DA, dark red), aggregated intraday position (IM, orange), imbalances (IB, black line), flexible demand consumption (dark green), battery behaviour (bright green), wind generation (light blue), PV generation (light green).}
    \label{FIG:EC_behaviour_summary}
\end{figure}

To honour the submitted bids, the energy community has to be able to supply or consume the committed electricity. The model provides the optimal behaviour of every element of the energy community to honour the bids at every scenario. Figure \ref{FIG:EC_behaviour_summary} shows a summary of the EC's market position at all markets considered and the optimal behaviour of all its elements in one scenario over a whole day.

\begin{figure}[ht]
 \begin{subfigure}{0.49\textwidth}
     \includegraphics[width=\textwidth]{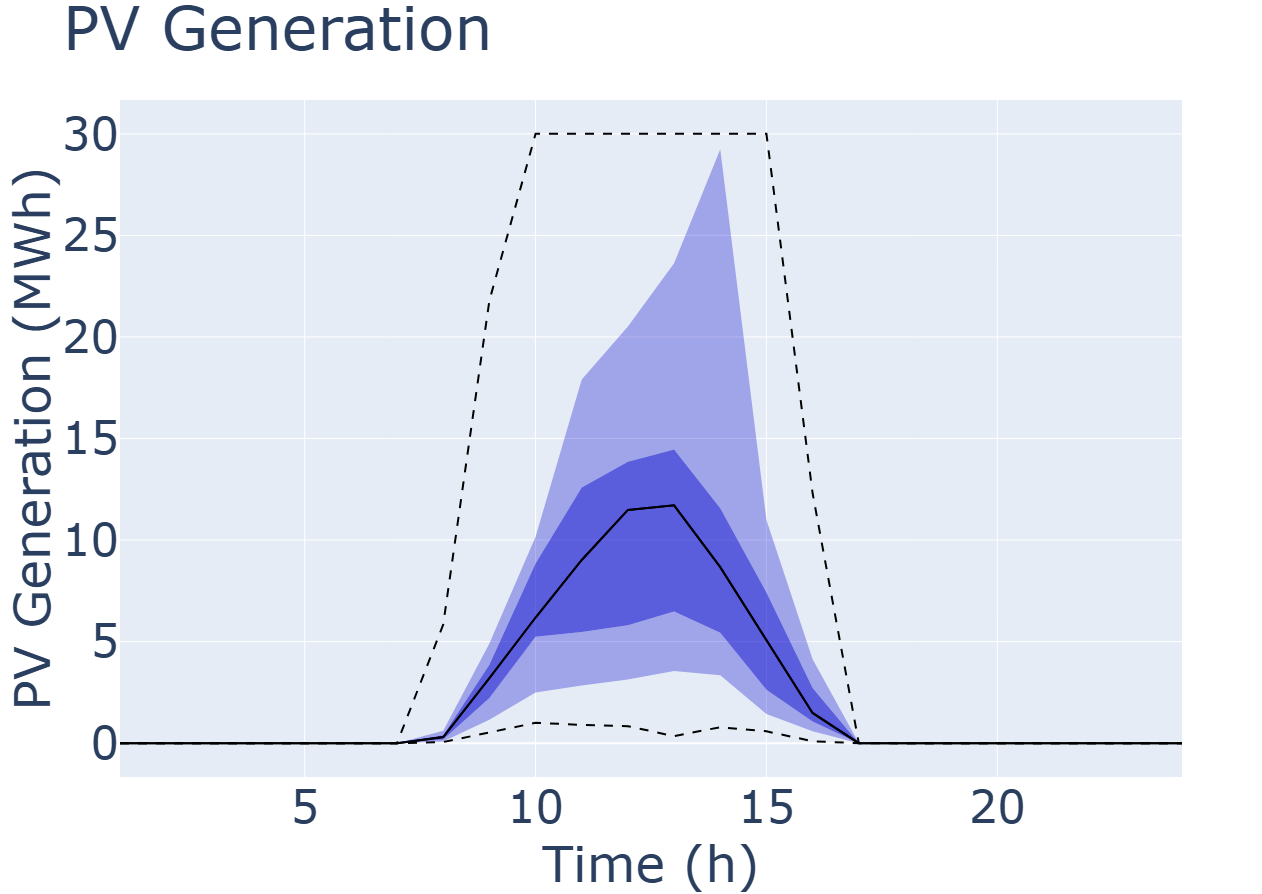}
 \end{subfigure}
 \hfill
 \begin{subfigure}{0.49\textwidth}
     \includegraphics[width=\textwidth]{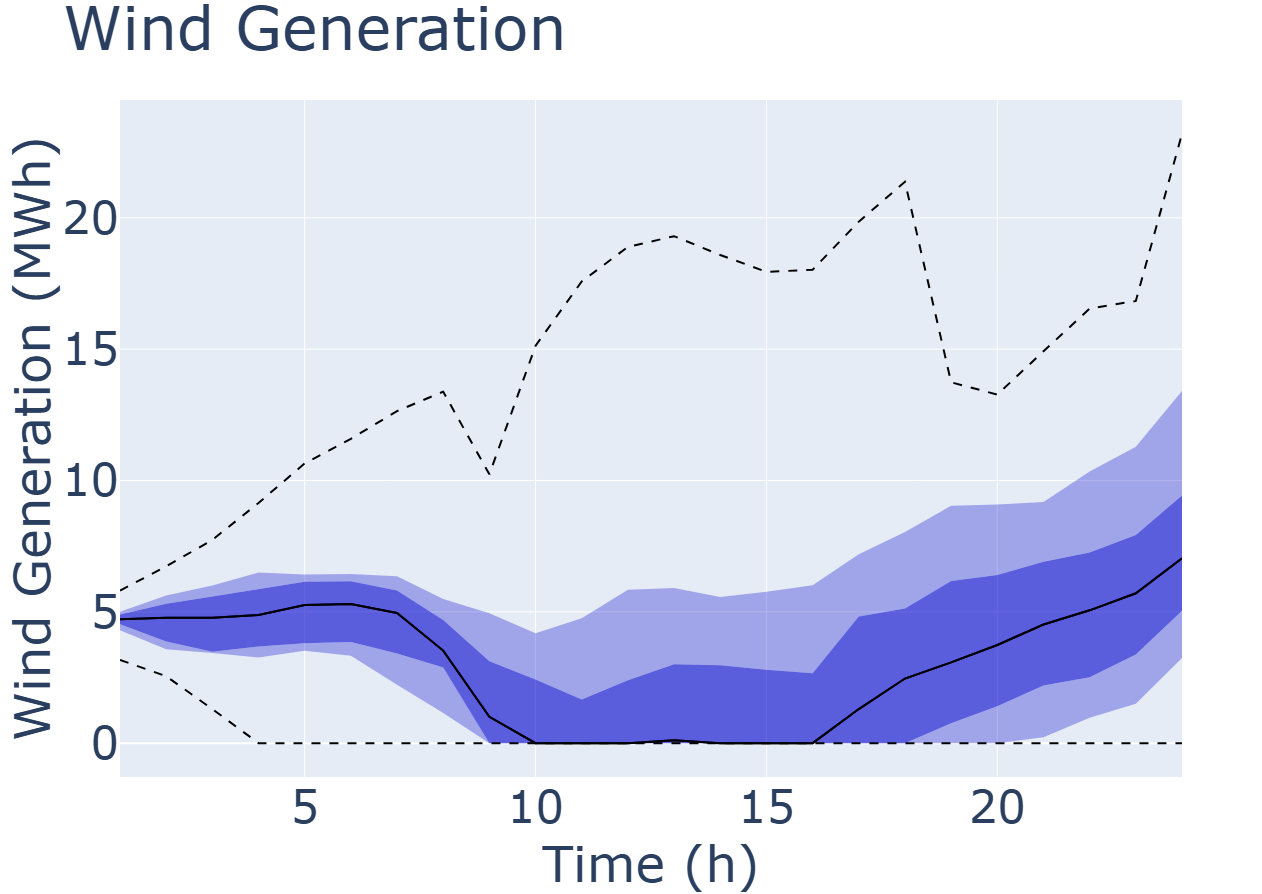}
 \end{subfigure}
 
 \medskip
 \begin{subfigure}{0.49\textwidth}
     \includegraphics[width=\textwidth]{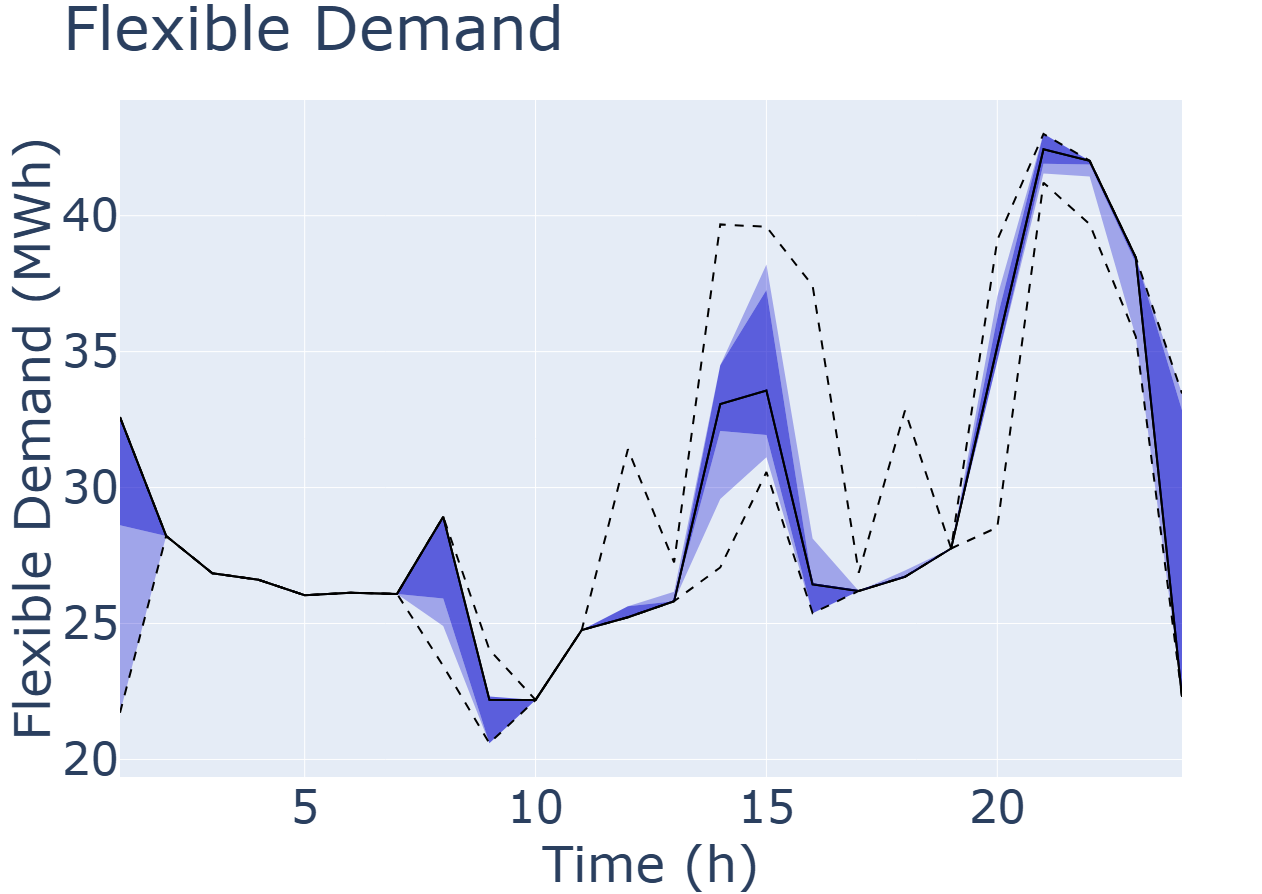}
 \end{subfigure}
 \hfill
 \begin{subfigure}{0.49\textwidth}
     \includegraphics[width=\textwidth]{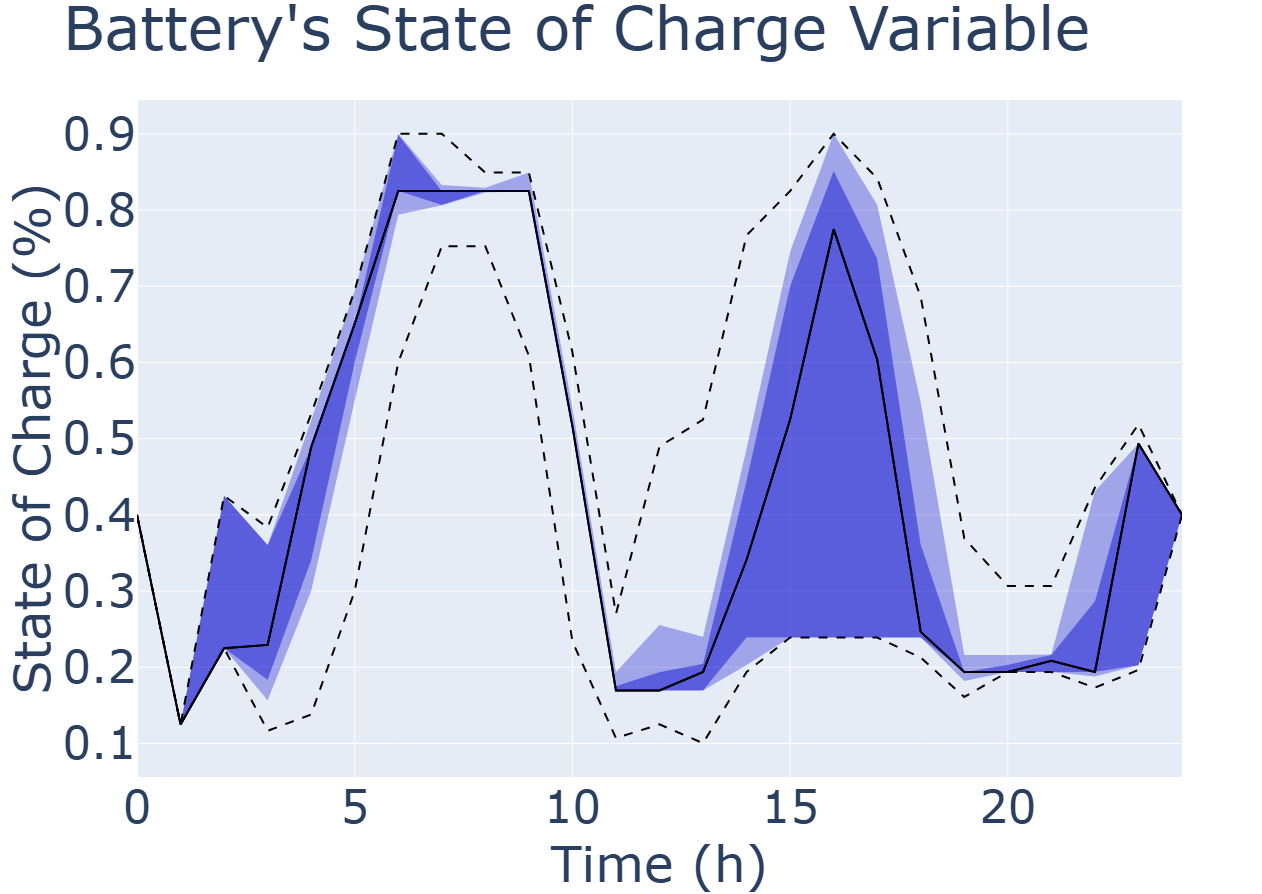}
 \end{subfigure}

 \caption{Percentiles of PV generation [MWh] (top left), wind generation [MWh] (top right), flexible electricity demand [MWh] (bottom left) and battery's state of charge [MWh] (bottom right) over a day. In every plot, there is: maximum and minimum (dashed black line), 10-90 percentile (light blue), 25-75 percentile (dark blue) and median (black line).}
 \label{FIG:EC_elements_behaviour}
\end{figure}

We can also analyse the behaviour of an energy community element across all scenarios. Figure \ref{FIG:EC_elements_behaviour} shows the 10, 25, 50, 75 and 90 percentiles of PV and wind generation, flexible demand, and the battery's state of charge across all scenarios in one day. The PV and wind generation data is obtained from the scenario generation procedure and are parameters of the optimization model, while the flexible demand position and the battery's state of charge is part of the optimal solution to the optimization problem.

\section{Analysis of Results} \label{SEC:AnalysisOfResults}

To analyse the model capabilities, we conduct a case study consisting on running the model for every day in December 2023 with real data. In this case study, the assets controlled by the energy community are solar PV panels with a total installed capacity of $\overline{P}^{PV} = 30$ MW, a wind farm of $\overline{P}^{W} = 30$ MW and a battery with an energy capacity of $E^{B} = 10$ MWh and maximum charging and discharging power rate of $P^{B} = 3$ MW.

The scenario trees of the model are built as described in Section \ref{SEC:GenerationOfScenariosTrees}. To fit the data, we use electricity price and variable renewable generation data from every hour in 2022 and 2023 up to the day for which the scenario is being generated. For example, to generate the scenario tree of the 16th December 2023, all hours between the 1st January 2022 and 15th December 2023 are used.

The data sources used to generate the scenarios are the following. Hourly electricity market prices for all markets considered have been obtained from OMIE, the Spanish and Portuguese market system operator \cite{OMIEData}. Solar PV generation is obtained from Renewables Ninja \cite{Pfenninger16, Staffell16, RenewablesNinja} in a location inside the Spanish postal code 08028, in Barcelona. Wind generation is also obtained from Renewables Ninja, but this time in a countryside location outside the city suitable for a small wind development. Moreover, the electricity demand data, which is not part of the scenario tree, is the hourly demand in December 2023 at the Spanish postal code 08028. It has been obtained from Datadis \cite{Datadis}.

The model outputs have been presented in Section \ref{SEC:ModelOutputs}. The figures presented there are from this case study. In the remaining of this section, we present aggregated results of this month-long experiment, and discuss their implications.

\begin{figure}[ht]
    \centering
    \includegraphics[width= \linewidth]{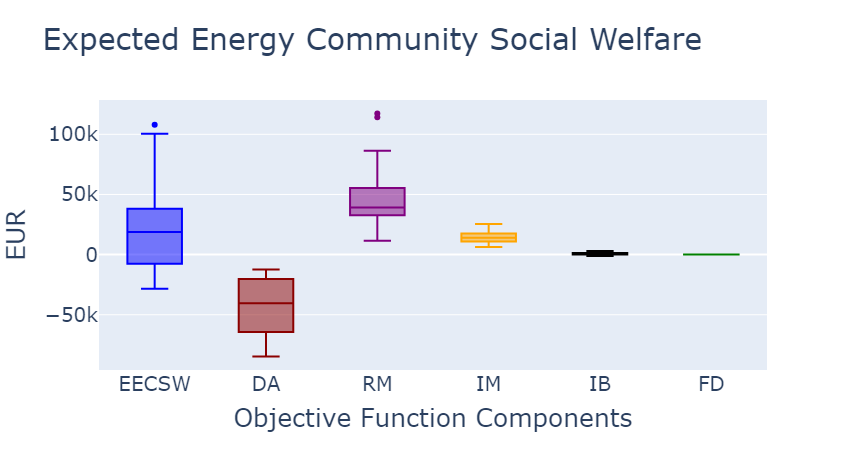}
    \caption{Expected Energy Community Social Welfare Analysis (EECSW) for the month of December 2023. Contribution of different objective function terms.}
    \label{FIG:EECSW_2024_12}
\end{figure}

Figure \ref{FIG:EECSW_2024_12} shows a box-plot of the objective function value across all the $31$ case study days, and the contribution of all its terms. EECSW is the objective function value, DA are the day-ahead incomes and costs \eqref{EQ:DAincome}, RM are the reserve market incomes \eqref{EQ:RESincome}, IM are all intraday market incomes and costs \eqref{EQ:IMincome}, IB are the imbalance collection rights and payment obligations \eqref{EQ:PosIBcollectRights}-\eqref{EQ:NegIBPayObligations}, and FD are the flexible demand penalties \eqref{EQ:DemFlexPenalty}.

Figure \ref{FIG:EECSW_2024_12} shows that the energy community has a positive social welfare about $75 \%$ of days. The largest contributor to this positive social welfare is the reserve market, highlighting the importance of energy community participation in this market as it can be a stream of large incomes. At the other end of the spectrum, the largest cost is incurred in the day-ahead market. This suggest that the energy community does not generate enough electricity to meet their daily electricity demand, but also that it might be willing to buy energy in the day-ahead market to make sure it can provide reserve at later stages. The fact that the net daily intraday market participation is always positive, i.e., always selling, supports this claim, as the energy community is willing to buy slightly more than it needs in the day-ahead market and recover part of this extra costs in the intraday markets.

Figure \ref{FIG:EECSW_2024_12} also shows that the imbalances incurred by the energy community are relatively small. The imbalance collection rights and payment obligations were also of the order of their sum, further supporting this claim. For this reason, they were not represented in Figure \ref{FIG:EECSW_2024_12} separately. Finally, the flexible demand also results in a small penalty, indicating that it is used in small quantities.

\section{Future Work}

In this section, we outline future research directions linked to this work. They can be classified in four main blocks: regulation adaptation and analysis, re-optimization at later stages, inclusion of the hydrogen chain, and open source release.

To better integrate European electricity markets, regulation is actively changing in Europe, and is expected to continue doing so in the near future. In particular, in recent years we have seen the redesign of balancing markets, the introduction of the continuous intraday electricity market and the move to $15$ minute periods in both balancing and intraday markets \cite{EU2019/943}. These was transposed to national legislations only recently, e.g., June 2024 in Spain \cite{BOE24_balances, BOE24_IM}. Moreover, day-ahead markets in Europe are expected to move to $15$ minute periods in late 2025 \cite{ENTSOE25}.

This greatly influences our model. First, the move to $15$ minute periods significantly increases the computational burden of the stochastic programming model because the number of decision variables indexed by time and the number of stages not linked to markets increase fourfold. The increase in the number of stages could be particularly critical for the model. Second, the continuous intraday electricity market clearing process is significantly different than the other electricity markets because it follows an order book format instead of an auction process. In particular, the uncertainties associated to this market need to be captured accurately.

Hence, in future work we will adapt the stochastic programming model to the new regulations with $15$ minute periods, and propose methodologies to address the challenges posed by the significant increase in problem size. Moreover, we will explore how to model the continuous intraday electricity market, and analyse the associated uncertainties.

A second future research direction is the re-optimization of the model at later stages. This is, to follow the decision process and re-optimize an adapted version of the model every time new information becomes available, i.e., every time a stage finishes. This will allow us to provide optimal actions for every decision the energy community makes with the most updated information available, and validate the methodology proposed with real data taking into account the complete decision process. Moreover, from a theoretical perspective, it is also interesting to analyse the value of the stochastic solution in a multi-stage stochastic programming environment.

A third future research direction is the inclusion of the hydrogen chain in the model. Hydrogen can be used as long-term energy storage, and can also contribute to the decarbonization of the heating and transport sector. Often, however, it is considered an expensive solution due to the low efficiencies in the generation of hydrogen from electricity and vice-versa. With the power of the presented market models, we will analyse the costs and benefits of including the hydrogen chain, allowing it to leverage revenue streams from multiple electricity markets.

Finally, a fourth future research direction is the implementation of the scenario generation methodology and multi-stage stochastic programming model in open source languages, and their release open source. In doing so, we will consider the inclusion of new features such as electricity demand uncertainty.

\section{Conclusions} \label{SEC:Conclusions}

In this work, we have presented a novel multi-stage stochastic programming model for the optimal participation of energy communities in the day-ahead, reserve and intraday electricity markets. The decision-making process is modelled in great detail, including a stage for every electricity market where the energy community might participate, and for every hour of the day, to accurately model the information available to make each decision and the hourly behaviour of the energy community. This results in a large multi-stage stochastic programming model with $34$ stages. 

To introduce the electricity price and variable renewable uncertainty into the model, scenario trees have been generated. The methodology to generate them has been adapted from previous work in the literature to fit the model and include hourly solar generation. Moreover, we have presented the three main model outputs in detail. These are the optimal bid of the energy community in the day-ahead market, the optimal price-accepting bid of the energy community in reserve and intraday markets, and the optimal behaviour of all elements of the energy community. Together, they provide the optimal day-ahead market bid of the energy community for every hour, and the optimal actions it should take at later stages to honour those bids.

A crucial aspect of an energy community is the simultaneous presence of an internal electricity demand, which can be flexible, and electricity generation capabilities. To capture this, we introduced a novel model for aggregated flexible demand at the community level. It is detailed enough to capture flexibility limits in the whole planning horizon and on specific periods, and has proven suitable to be incorporated in the complex stochastic programming problem considered. 

Furthermore, the model allows the energy community to provide buying, selling and combined buying and selling bids. The combined buying and selling bids are bids where the energy community buys if prices are below a certain threshold and sells if they are above that threshold, which, to the best of our knowledge, have attracted little attention in the literature. In our analysis, we observed that the optimal solution contained all bid types, depending on the day and time period. This highlights the importance of giving the energy community the flexibility to choose which type of bid to submit to maximize its social welfare.

The case study presented, with data from December 2023, shows that the energy community is able to achieve a positive social welfare in almost $75\%$ of days. Interestingly, in the days considered, the energy community always incurs a cost in the day-ahead market while the reserve market is the most important source of incomes. This suggests that the energy community is willing to buy more energy than it needs in the day-ahead market to make sure it can provide reserve, which seems to be more lucrative. It also underscores the importance of the reserve market as a source of income.

Finally, highlight that, while the model presented here has been customized for energy communities, it can be adapted to any entity that has an internal electricity demand and variable renewable generation and is able and willing to participate in electricity markets. Examples of such entities could be aggregators or owners of large electrical energy storage systems such as large batteries or pumped hydro storage.

\section*{Acknowledgements}

We thank Andrea Ademollo, Jan Jettmann and Cristian Pachón for their valuable contributions, discussions and insights in the preparation of this work.

\subsection*{CRediT Authorship Contribution Statement}

\textbf{Albert Sol\`a Vilalta}: Formal analysis, Methodology, Software, Writing – original draft, Writing - review and editing. \textbf{Ignasi Ma\~{n}\'e Bosch:} Formal analysis, Methodology, Software. \textbf{F.- Javier Heredia:} Formal analysis, Funding acquisition, Methodology, Supervision, Writing – review and editing. 

\subsection*{Funding Sources}

This work has been supported with grants TED2021-131365B-C44 and Cetp-FP-2023-00185 from the Spanish Ministerio de Ciencia, Innovación y Universidades.

\appendix
\section{Mathematical Formulation of the Multi-Stage Stochastic Programming Model}
\label{APP:MathematicalFormulation}

In this appendix, we present the mathematical formulation of the multi-stage stochastic programming model introduced in this work. Since it is quite a large model, the sets, parameters, variables, and constraints are defined on separate sections.

We consider an energy community with the following elements:
\begin{enumerate}
	\item (Implicit) Flexible Demand.
	\item Wind Farm.
	\item Solar Farm.
	\item Battery Energy Storage System (BESS).
\end{enumerate}

\subsection{Fundamental Elements} \label{SUBSEC:FundamentalElements}
\def\Ts{\mathcal{T}}
\def\Is{\mathcal{I}}
\def\Os{\Omega}
\def\Ss{\mathcal{S}}
\def\Cs{\mathcal{C}}
\begin{tabularx}{\textwidth}{l l}
    \textbf{Sets:} \\
    $\mathcal{T}$       & set of time periods [hours]. \\
	$\mathcal{T}_0$     & extended ordered set of time periods $\mathcal{T}_0 = \{ 0\}\cup\mathcal{T}$ [hours]. \\     
    $\mathcal{I}$       & set of intraday markets. \\    
    $\mathcal{I}_t \subseteq \mathcal{I}$
                        & set of intraday markets where $t \in \mathcal{T}$ is a bidding period. \\   
    $\mathcal{T}_i \subseteq \mathcal{T}$
                        & set of bidding periods in intraday market $i \in \mathcal{I}$ [hours]. \\
    $\Omega$            & set of scenarios. \\
    $\mathcal{S}$       & set of stages. \\
	$\mathcal{S}_0$     & extended set of stages, $\mathcal{S}_0=\{0\}\cup\mathcal{S}$. \\
    $\Cs_s$             & set of clusters at stage $s\in\Ss$.\\
    $\Cs_{s,c}$         & set of scenarios of cluster $c\in \Cs_s$ at stage $s\in\Ss$, \\
    &                   $\Cs_{s,c}=\{\omega_1,\omega_2,\ldots,\omega_{|\Cs_{s,c}|}\}$.\\
    \\
\end{tabularx}
\begin{tabularx}{\textwidth}{l l}
    \textbf{Parameters:} \\
	$S^{IM}_i$ & stage where the price of intraday market $i \in \Is$ is revealed. \\
	$\Pi_{\omega}$& probability of scenario $\omega \in \Omega$. \\
    \\
\end{tabularx}


\subsection{Modelling of an Energy Community} 
\label{SUBSEC:TheEnergyCommunityModel}

\subsubsection{Flexible Demand} 
\label{SUBSEC:FlexibleDemand}

\begin{tabularx}{\textwidth}{l l}
    \textbf{Sets:} \\
    $\mathcal{F}$& set of time period intervals where a given flexible demand \\
    & has to be met. \\ 
    \\
\end{tabularx}
\begin{tabularx}{\textwidth}{l l}
    \textbf{Parameters:} \\

    $D_t$& central demand at $t \in \mathcal{T}$ [MWh]. \\
    $\underline{D}_t$& minimum demand at $t \in \mathcal{T}$ (inflexible demand) [MWh]. \\
    $\overline{D}_t$& maximum demand at $t \in \mathcal{T}$ [MWh]. \\
    $T_f^-$& first time step of interval $f \in \mathcal{F}$ [hour]. \\
    $T_f^+$& last time step of interval $f \in \mathcal{F}$ [hour]. \\
    $c_f$ & fraction of central demand that must be met in the \\
    & interval $f \in \mathcal{F}$ [MWh]. \\
    $C^{FD}$& cost per unit of flexible demand used [\euro/MWh]. \\
    \\
\end{tabularx}
\begin{tabularx}{\textwidth}{l l}
    \textbf{Variables:} \\
    $f_{t, \omega}$& flexible demand at $t \in \mathcal{T}$ under scenario $\omega \in \Omega$ (after \\ 
    & flexibility is considered) [MWh]. \\
    $f_{t, \omega}^+$& auxiliary variable. Positive part of flexible demand shift \\ 
    & from central demand $D_t$ [MWh]. \\
    $f_{t, \omega}^-$& auxiliary variable. Negative part of flexible demand shift \\
    & from central demand $D_t$ [MWh]. \\
    \\
    \end{tabularx}

First, we introduce nonnegativity constraints to the flexible demand variables
\begin{equation} \label{EQ:FlexDemandPos}
	f_{t, \omega} \geq 0, \ f_{t, \omega}^+ \geq 0, \ f_{t, \omega}^- \geq 0 \quad \forall t \in \mathcal{T}, \ \forall \omega \in \Omega.
\end{equation}

The flexible demand is bounded by the minimum demand and maximum demand at every time step. This is
\begin{equation} \label{EQ:FlexDemand_UB_preliminary}
	f_{t, \omega} \leq \overline{D}_t \quad \forall t \in \mathcal{T}, \ \forall \omega \in \Omega
\end{equation}
and 
\begin{equation} \label{EQ:FlexDemand_LB_preliminary}
	\underline{D}_t \leq f_{t, \omega} \quad \forall t \in \mathcal{T}, \ \forall \omega \in \Omega.
\end{equation}
These two constraints are weaker than their counterparts \eqref{EQ:FlexDemand_UB} and \eqref{EQ:FlexDemand_LB}. \eqref{EQ:FlexDemand_UB} and \eqref{EQ:FlexDemand_LB} guarantee that, if reserve has to be activated, flexible demand's behaviour plus reserve commitments still result in feasible actions at any time period. Hence, we do not introduce \eqref{EQ:FlexDemand_UB_preliminary} nor \eqref{EQ:FlexDemand_LB_preliminary} in the model's implementation, but only describe them here to present the model.

Moreover, the total daily demand has to meet the central estimate of the demand over the day: 
\begin{equation} \label{EQ:DailyDem}
	\sum_{t \in \mathcal{T}} f_{t, \omega} = \sum_{t \in \mathcal{T}} D_t \quad \forall \omega \in \Omega.
\end{equation}
This is done to avoid meeting a lower demand over the full time horizon of the problem ($24$h), or a higher demand if the reward for downward reserve is consistently high and/or electricity prices are consistently negative. We also introduce time intervals where a certain percentage of the central demand needs to be met over that interval:
\begin{equation} \label{EQ:InterDem}
	\sum_{T_f^- \leq t \leq T_f^+} f_{t, \omega} \geq c_f \sum_{T_f^- \leq t \leq T_f^+} D_t \quad \forall f \in \mathcal{F}, \ \forall \omega \in \Omega.
\end{equation}

To model the cost of implicit flexible demand, we introduce a penalty in the objective function to use it. This requires to split the variables $f_{t, \omega}$ into its positive and negative parts, which is done in the following constraints
\begin{equation} \label{EQ:FlexDemDisplace}
	D_t - f_{t, \omega} = f_{t, \omega}^+ - f_{t, \omega}^- \quad \forall t \in \mathcal{T}, \ \forall \omega \in \Omega.
\end{equation}

\subsubsection{Renewable Generation: PV and Wind Farm} 
\label{SUBSUBSECTION:SolarFarm}

The renewable production of the energy community, solar photovoltaic and wind, are considered in this model as random parameters with values between zero and some given maximum production $\bar{P}^{PV}$ and $\bar{P}^{PV}$ respectively.

\begin{tabularx}{\textwidth}{l l}
    \textbf{Parameters:} \\
	$PV_{t, \omega}$       & solar PV generation at $t \in \mathcal{T}$ under scenario $\omega \in \Omega$ [MWh]. \\
	$\overline{P}^{PV}$    & solar PV nameplate capacity [MW]. \\
	$W_{t, \omega}$        & wind generation at time $t \in \mathcal{T}$ under scenario $\omega \in \Omega$ [MWh]. \\
    $\overline{P}^W$       & wind nameplate capacity [MW]. \\
    $S^{R}_t$              & stage where renewable production (solar and wind) at time \\ 
                           & $t\in\Ts$ is observed. \\  
	\\  
\end{tabularx}

We assume that observations of the solar and wind random parameters at time $t \in \mathcal{T}$ are made at the same stage $S^{R}_t$.

\subsubsection{Battery Energy Storage System (BESS)}
\label{SUBSUBSEC:BatteryEnergyStorageSystem}

\begin{tabularx}{\textwidth}{l l}
    \textbf{Parameters:} \\
	$E^B$              & BESS energy capacity [MWh].	 \\
	$P^B$              & BESS maximum charging/discharging rate [MW]. \\
	$\eta^B$           & BESS round-trip efficiency. \\
	$\overline{\SOC}$  & BESS maximum state of charge. Unitless, $\overline{\SOC} \in [0, 1]$. \\
	$\underline{\SOC}$ & BESS minimum state of charge. Unitless, $\underline{\SOC} \in [0, 1]$. \\
	$\SOC_0$           & BESS initial state of charge. Unitless, $\SOC_0 \in [0, 1]$. \\
	$\SOC_T$           & BESS final state of charge. Unitless, $\SOC_T \in [0, 1]$. \\
	\\
       
    \textbf{Variables:} \\
	$c_{t, \omega}$    & BESS charge rate at time $t \in \mathcal{T}$ under scenario \\
                       & $\omega \in \Omega$ [MW]. \\
	$d_{t, \omega}$    & BESS discharge rate at time $t \in \mathcal{T}$ under scenario \\
                       & $\omega \in \Omega$ [MW]. \\	
	$id_{t, \omega}$   & binary. $id_{t, \omega} = 0$ (resp. $id_{t, \omega}= 1$) if BESS is charging  \\
	&  (resp. discharging) at time $t \in \mathcal{T}$ under scenario $\omega \in \Omega$. \\
	$\soc_{t, \omega}$ & state of charge of BESS at time $t \in \mathcal{T}_0$ under scenario \\
	& $\omega \in \Omega$. Unitless, $\soc_{t, \omega} \in [0, 1]$. \\
	\\
\end{tabularx}

First, we introduce the following nonnegativity  constraints, and define $id_{t, \omega}$ as a binary:
\begin{equation} \label{dVdCPos}
	c_{t, \omega} \geq 0, \quad d_{t, \omega} \geq 0, \quad id_{t, \omega} \in \{ 0, 1 \} \quad \forall t \in \mathcal{T}, \ \forall \omega \in \Omega.
\end{equation}

We introduce limits in the charged and discharged quantities of the battery, as well as not allowing simultaneous charging and discharging in every time step:
\begin{equation} \label{EQ:VPP_state_dV}
	d_{t, \omega} \leq P^B id_{t, \omega} \quad \forall t \in \mathcal{T}, \ \forall \omega \in \Omega
\end{equation}
and
\begin{equation} \label{EQ:VPP_state_cV}
	c_{t, \omega} \leq P^B (1 - id_{t, \omega}) \quad \forall t \in \mathcal{T}, \ \forall \omega \in \Omega.
\end{equation}

We define the variable state of charge $\soc_{t, s}$ of the BESS with constraints: 
\begin{equation} \label{EQ:SOCV}
	\soc_{t, \omega} = \soc_{t - 1, \omega} + \frac{c_{t, \omega} - \frac{d_{t, \omega}}{\eta^B}}{E^B} \quad \forall t \in \mathcal{T}, \ \forall \omega \in \Omega. 
\end{equation} 
They have to be within the operational limits of the BESS:
\begin{equation} \label{EQ:socVOperLimits}
	\underline{\SOC} \leq \soc_{t, \omega} \leq \overline{\SOC} \quad \forall t \in \mathcal{T}, \ \forall \omega \in \Omega.
\end{equation}
The initial and final states of charge of the BESS are defined with:
\begin{equation} \label{EQ:SOCV_ini}
	\soc_{0, \omega} = \SOC_0 \quad \forall \omega \in \Omega
\end{equation}
and
\begin{equation} \label{EQ:SOCV_fin}
	\soc_{T, \omega} = \SOC_T \quad \forall \omega \in \Omega.
\end{equation}

\subsection{Modelling Electricity Market Participation} 
\label{SUBSEC:ModellingElectricityMar}

\subsubsection{Day-Ahead Market}
\label{SUBSUBSEC:DayAheadMarket}
\def\ODAs{\Omega^{DA}}
\def\oda{\omega^{DA}}

\begin{tabularx}{\textwidth}{l l}
    \textbf{Sets:} \\
    $\ODAs$ & Scenarios in $j\in\Os$ sorted in ascending order of day-ahead prices,\\
            &$\ODAs=\{\omega_1,\omega_2,\ldots,\omega_{|\Os|}\}$, $\lambda^{DA}_{t,\omega_i}<\lambda^{DA}_{t,\omega_{i+1}}$, $\forall\omega_{i<|\Os|}$.\\
	& For every $t \in \mathcal{T}$, provides the scenario pairs $(l, j)$ s.t. $\lambda^D_{t, l} \leq \lambda^D_{t, j}$. \\	
	& Used to ensure monotonicity of DA bids. \\
\end{tabularx}

\begin{tabularx}{\textwidth}{l l}
	\\
    \textbf{Parameters:} \\
	$\lambda^{D}_{t, \omega}$& DA market price at $t \in \mathcal{T}$ under scenario $\omega \in \Omega$ [\euro/MWh]. \\
    $\underline{P}^{DA}$ & minimum bid size to DA market [MWh]. \\
    \\
\end{tabularx}

\begin{tabularx}{\textwidth}{l l}
    \textbf{Variables:} \\
	$e^{DA+}_{t, \omega}$& DA energy matched of an EC selling bid at time $t \in \mathcal{T}$ \\
	& under scenario $\omega \in \Omega$ [MWh]. \\
	$e^{DA-}_{t, \omega}$ & DA energy matched of an EC buying bid at time $t \in \mathcal{T}$ \\
	& under scenario $\omega \in \Omega$ [MWh]. \\
	$ie^{DA+}_{t, \omega}$& binary. $ie^{DA+} = 1$ if the EC sells electricity at time $t \in \mathcal{T}$  \\ 
    &  under scenario $\omega \in \Omega$. $ie^{DA+}_{t, \omega} = 0$ o/w. \\
    $ie^{DA-}_{t, \omega}$ & binary. $ie^{DA-}_{t, \omega} = 1$ if the EC buys electricity at time $t \in \mathcal{T}$ \\
	& under scenario $\omega \in \Omega$. $ie^{DA-}_{t, \omega} = 0$ o/w.
    \\
\end{tabularx}

First, we impose that the variables $e^{DA+}_{t, \omega}, e^{DA-}_{t, \omega}$ are nonegative, and that $ie^{DA+}_t, ie^{DA-}_t$ are binary:
\begin{equation} \label{EQ:DAVarDef}
	e^{DA+}_{t, \omega} \geq 0, \ e^{DA-}_{t, \omega} \geq 0, \ ie^{DA+}_t, ie^{DA-}_t \in \{0, 1 \} \quad \forall t \in \mathcal{T} \quad \forall \omega \in \Omega.
\end{equation}

We introduce lower and upper bounds on the selling and buying energy matched quantities represented by the variables $e^{DA+}_{t, \omega}$ and $e^{DA-}_{t, \omega}$. Moreover, at every scenario $\omega \in \Omega$, the market position of the EC has to be to sell, buy or neither of the two. This is guaranteed with the binary variables $ie^{DA+}_{t, \omega}$ and $ie^{DA-}_{t, \omega}$. The lower bounds of \eqref{EQ:DA_sell_bid_LB} and \eqref{EQ:DA_buy_bid_LB} represent the minimum bid size. The upper bounds of \eqref{EQ:DA_sell_bid_LB} and \eqref{EQ:DA_buy_bid_LB} are such that all elements of the energy community can contribute towards the day-ahead market bid. Note also that speculative behaviour is controlled because quantity bids cannot be larger than the technical capacity of the energy community. For selling energy matched, we have 
\begin{equation} \label{EQ:DA_sell_bid_LB}
    \underline{P}^{DA} ie^{DA+}_{t, \omega} \leq e^{DA+}_{t, \omega} \leq (\overline{P}^W + \overline{P}^{PV} + P^B - \underline{D}_t) ie^{DA+}_{t, \omega}  \quad \forall t \in \mathcal{T}, \quad \forall \omega \in \Omega.
\end{equation}
and for buying energy matched 
\begin{equation} \label{EQ:DA_buy_bid_LB}
	\underline{P}^{DA} ie^{DA-}_{t, \omega} \leq e^{DA-}_{t, \omega} \leq (P^B + \overline{D}_t) ie^{DA-}_{t, \omega} \quad \forall t \in \mathcal{T}, \quad \forall \omega \in \Omega.
\end{equation}

At every scenario, we also impose that it is not possible to simultaneously have matched bought and sold energy quantities (\verb|DA_buy_or_sell|)
\begin{equation} \label{EQ:DA_buy_or_sell}
	ie^{DA+}_{t, \omega} + ie^{DA-}_{t, \omega} \leq 1 \quad \forall t \in \mathcal{T}, \quad \forall \omega \in \Omega.
\end{equation}
Note that this allows the EC to submit neither of the two.

Finally, to be able to construct a bidding curve for the EC from the results obtained at every scenario, the quantities offered in electricity markets have to be monotone w.r.t. the scenario prices. For the day-ahead market, this is guaranteed with the following constraints 
\begin{equation} \label{EQ:DA_bid_mono_1}
    e^{DA+}_{t,\oda_i}\le e^{DA+}_{t,\oda_{(i+1)}} \quad \forall  \oda_{i<|\Os|}\in\ODAs, t\in\Ts
\end{equation}
\begin{equation} \label{EQ:DA_bid_mono_2}
    e^{DA-}_{t,\oda_i}\le e^{DA-}_{t,\oda_{(i+1)}} \quad \forall  \oda_{i<|\Os|}\in\ODAs, t\in\Ts
\end{equation}

More specifically, in scenarios where the EC sells electricity, the higher the electricity price gets, the larger the selling electricity matched should be. Symmetrically, in scenarios where the EC buys electricity, the higher the electricity price gets, the smaller the buying electricity matched should be. 

\subsubsection{Reserve Market} 
\label{SUBSUBSEC:ReserveMarket}

We assume that the wind and solar farms cannot contribute to provide reserve. The other elements of the energy community (Battery and Flexible Demand) can. This follows the approach by \cite{Heredia18}, where VPP's bid to reserve markets is limited by the battery capabilities only.

We model the participation of the EC in the reserve market with a price-accepting bid for every day-ahead scenario cluster. There are different price-accepting bids because the day-ahead market prices have already been observed at the time of making the reserve market decisions. Note that this is different to the modelling of day-ahead bids presented in Section \ref{SUBSUBSEC:DayAheadMarket}, where a complete day-ahead bid curve is obtained and it is unique across all scenarios.

\begin{tabularx}{\textwidth}{l l}
    \textbf{Parameters:} \\
    	$\lambda^{R}_{t, \omega}$ & RM price at time $t \in \mathcal{T}$ under scenario $\omega \in \Omega$ [\euro/MW]. \\
    	$T^R$& duration of reserve response [h]. \\
    	$R^{U, FD}_t$& upper bound on the upwards reserve provided by \\ 
        & flexible demand [MW]. \\
    	$R^{D, FD}_t$& upper bound on the downwards reserve provided by \\
        & flexible demand [MW]. \\
    \\
\end{tabularx}

\begin{tabularx}{\textwidth}{l l} 
    \textbf{Variables:} \\
	$r^{U}_{t, \omega}$& upwards reserve matched by the EC at time $t \in \mathcal{T}$ under \\
    & scenario $\omega \in \Omega$ [MW]. \\
    $r^{D}_{t, \omega}$& downwards reserve matched by the EC at time $t \in \mathcal{T}$ under \\
    & scenario $\omega \in \Omega$ [MW]. \\
    $r^{U, B}_{t, \omega}$& BESS component of upward reserve matched at time $t \in \mathcal{T}$ \\
    & under scenario $\omega \in \Omega$ [MW]. \\
    $r^{D, B}_{t, \omega}$ & BESS component of downward reserve matched at time \\
    & $t \in \mathcal{T}$ under scenario $\omega \in \Omega$ [MW]. \\
    $r^{U, FD}_{t, \omega}$& flexible demand component of upward reserve matched \\ 
    & at time $t \in \mathcal{T}$ under scenario $\omega \in \Omega$ [MW]. \\
    $r^{D, FD}_{t, \omega}$& flexible demand component of downward reserve matched \\
    & at time $t \in \mathcal{T}$ under scenario $\omega \in \Omega$ [MW]. \\
    \\
\end{tabularx}

First, we impose that the reserve variables are nonnegative
\begin{equation} \label{EQ:RVarDef}
    r^{U}_{t, \omega} \geq 0, \quad r^{D}_{t, \omega} \geq 0, \quad r^{U, B}_{t, \omega} \geq 0, \quad r^{D, B}_{t, \omega} \geq 0, \quad r^{U, FD}_{t, \omega} \geq 0, \quad r^{D, FD}_{t, \omega} \geq 0
\end{equation}
for all $t \in \mathcal{T}$, and $\omega \in \Omega$. 

The reserve provided by the EC is the sum of the reserve provided by those elements able to provide reserve. This is
\begin{equation} \label{EQ:RUcomp}
    r^{U}_{t, \omega} = r^{U, B}_{t, \omega} +  r^{U, FD}_{t, \omega} \quad \forall t \in \mathcal{T}, \ \forall \omega \in \Omega
\end{equation}
and
\begin{equation} \label{EQ:RDcomp}
	r^{D}_{t, \omega} = r^{D, B}_{t, \omega} +  r^{D, FD}_{t, \omega} \quad \forall t \in \mathcal{T}, \ \forall \omega \in \Omega.
\end{equation}

We ensure that the elements of the energy community capable of providing reserve would be capable to satisfy their reserve commitments if instructed to do so. 

For flexible demand, we ensure that there is enough flexibility (in energy terms) at every time step to satisfy the reserve requirements through
\begin{equation} \label{EQ:FlexDemand_UB}
	f_{t, \omega} + T^R r^{D, FD}_{t, \omega} \leq \overline{D}_t \quad \forall t \in \mathcal{T}, \ \forall \omega \in \Omega
\end{equation}
and
\begin{equation} \label{EQ:FlexDemand_LB}
	\underline{D}_t \leq f_{t, \omega} - T^R r^{U, FD}_{t, \omega} \quad \forall t \in \mathcal{T}, \ \forall \omega \in \Omega.
\end{equation}

We also ensure that there is enough flexibility (in power terms) at every time step to satisfy the reserve requirements using
\begin{equation} \label{EQ:FlexDemP_DReserve}
	r^{D, FD}_{t, \omega} \leq R^{D, FD}_t \quad \forall t \in \mathcal{T}, \ \forall \omega \in \Omega
\end{equation}
and
\begin{equation} \label{EQ:FlexDemP_UReserve}
	r^{U, FD}_{t, \omega} \leq R^{U, FD}_t \quad \forall t \in \mathcal{T}, \ \forall \omega \in \Omega.
\end{equation}

The operation of the battery also needs to take into account the committed reserve. For the charging and discharging rates, we have:
\begin{equation} \label{EQ:VPP_RU2s}
	r^{U, B} _{t, \omega} - c_{t, \omega} + d_{t, \omega}  \leq P^B \quad \forall t \in \mathcal{T}, \ \forall \omega \in \Omega
\end{equation}
and
\begin{equation} \label{EQ:VPP_RD2s}
	r^{D, B}_{t, \omega} + c_{t, \omega} - d_{t, \omega} \leq P^B \quad \forall t \in \mathcal{T}, \ \forall \omega \in \Omega.
\end{equation}

The reserve provided by the battery also requires a sufficiently high (resp. low) state of charge so that the upwards (resp. downards) reserve can be provided, if needed. The following constraints ensure it for one time step:
\begin{equation} \label{EQ:VPP_RU_SOCV}
	\underline{\SOC} \leq \soc_{t, \omega} - \frac{T^R r^{U, B}_{t, \omega}}{\eta^B E^B} \quad \forall t \in \mathcal{T}, \ \forall \omega \in \Omega
\end{equation}
and
\begin{equation} \label{EQ:VPP_RD_SOCV}
	\soc_{t, \omega} + \frac{T^R r^{D, B}_{t, \omega} }{E^B} \leq \overline{\SOC} \quad \forall t \in \mathcal{T}, \ \forall \omega \in \Omega.
\end{equation}

\subsubsection{Intraday Markets}
\label{SUBSUBSEC:IntradayMarkets}

We model the participation of the EC in intraday markets with price-accepting bids. There is a different price-accepting bid for every scenario cluster of the stage immediately before the corresponding intraday market stage. Note that this is similar to the modelling of reserve market bids in Section \ref{SUBSUBSEC:ReserveMarket}, with the corresponding stage adaptation, and different to the modelling of complete day-ahead market bids presented in Section \ref{SUBSUBSEC:DayAheadMarket}.

\begin{tabularx}{\textwidth}{l l}    
    \textbf{Parameters:} \\
    	$\lambda^I_{i, t, \omega}$& intraday $i \in \mathcal{I}$ price at time $t \in \mathcal{T}$ under scenario \\
        & $\omega \in \Omega$ [\euro/MWh]. \\
	$R^{IDA}$ & bound on the ratio of energy traded in intraday \\
	& markets and day-ahead market (to control speculation).\\
    \\ 
\end{tabularx}

\begin{tabularx}{\textwidth}{l l}     
    \textbf{Variables:} \\
    $e^{IM}_{i, t, \omega}$& energy matched by the EC at intraday market $i \in \mathcal{I}$ \\
    & at time $t \in \mathcal{T}$ under scenario $\omega \in \Omega$ [MWh]. \\
	\\
\end{tabularx}

To avoid speculative behaviour, the bid of the EC to intraday markets should be small compared to the day-ahead market
\begin{equation} \label{EQ:DM_IM_1}
	-R^{IDA} (e^{DA+}_{t, \omega} + e^{DA-}_{t, \omega})  \leq \sum_{i \in \mathcal{I}_t} e^{IM}_{i, t, \omega} \leq R^{IDA} (e^{DA+}_{t, \omega} + e^{DA-}_{t, \omega}) \quad \forall t \in \mathcal{T}, \ \forall \omega \in \Omega
\end{equation}
and
\begin{equation} \label{EQ:IM_bounds_3-4}
	-R^{IDA} (e^{DA+}_{t, \omega} + e^{DA-}_{t, \omega}) \leq e^{IM}_{i, t, \omega} \leq R^{IDA} (e^{DA+}_{t, \omega} + e^{DA-}_{t, \omega}) \quad \forall i \in \mathcal{I}, \ \forall t \in \mathcal{T}_i, \ \forall \omega \in \Omega.
\end{equation}


\subsubsection{System Imbalances}
\label{SUBSUBSEC:SystemImbalances}

\begin{tabularx}{\textwidth}{l l}
    \textbf{Parameters:} \\
    $\lambda^{IB+}_{t, \omega}$& positive imbalance price at time $t \in \mathcal{T}$ under scenario \\
    & $\omega \in \Omega$ [\euro/MWh]. \\
	$\lambda^{IB-}_{t, \omega}$ & negative imbalance price at time $t \in \mathcal{T}$ under scenario \\
    & $\omega \in \Omega$ [\euro/MWh]. \\
    $E^{IB+}_{t, \omega}$& upper bound on the positive imbalance at time $t \in \mathcal{T}$ \\
    & under scenario $\omega \in \Omega$ [MWh]. \\
	$E^{IB-}_{t, \omega}$& upper bound on the negative imbalance at time $t \in \mathcal{T}$ \\ 
	& under scenario $\omega \in \Omega$ [MWh]. \\
    \\ 
       
    \textbf{Variables:} \\
    $e_{t, \omega}^{IB-}$& negative imbalance of EC at time $t \in \mathcal{T}$ under scenario \\
    & $\omega \in \Omega$ [MWh].	\\
	$e_{t, \omega}^{IB+}$& positive imbalance of EC at time $t \in \mathcal{T}$ under scenario \\
    & $\omega \in \Omega$ [MWh].	\\ 
    \\
\end{tabularx}

First, we impose that the imbalance variables are non-negative  
\begin{equation} \label{EQ:IBVarDef}
    e_{t, \omega}^{IB+} \geq 0, \quad e_{t, \omega}^{IB-} \geq 0 \quad \forall t \in \mathcal{T}, \ \forall \omega \in \Omega.
\end{equation}

The imbalances of the EC are defined in the following equation
\begin{equation} \label{EQ:IBDef}
    e^{IB+}_{t, \omega} - e^{IB-}_{t, \omega} = e^{DA-}_{t, \omega} + W_{t, \omega} + PV_{t, \omega} + d_{t, \omega} - \Big (e^{DA+}_{t, \omega} + \sum_{i \in \mathcal{I}_t} e^{IM}_{i, t, \omega} + f_{t, \omega} + c_{t, \omega} \Big ) \quad \forall t \in \mathcal{T}, \ \forall \omega \in \Omega.
\end{equation}

We introduce upper and lower bounds to the imbalances, to avoid speculative behaviour:
\begin{equation} \label{EQ:VPP_IB_MAX}
	e_{t, \omega}^{IB+} \leq E_{t, \omega}^{IB+} \quad \forall t \in \mathcal{T}, \ \forall \omega \in \Omega
\end{equation}
and
\begin{equation} \label{EQ:VPP_IB_MIN}
	e_{t, \omega}^{IB-} \leq E_{t, \omega}^{IB-} \quad \forall t \in \mathcal{T}, \ \forall \omega \in \Omega.
\end{equation}

\subsection{The Nonanticipativity Constraints} \label{SUBSEC:TheNonAnticipativityConstraints}

\label{SEC:NonanticipativityConstraints}

The nonanticipativity constraints ensure that every decision is made with the amount of information available at the time of making that decision.

For the day-ahead and reserve markets, the nonanticipativity constraints are:
\begin{subequations}
    \begin{empheq}
    [right=\empheqrbrace
    {,\ \omega_{i<|\mathcal{C}_{1,c}|}\in\mathcal{C}_{1,c}, c \in \mathcal{C}_1, t\in\mathcal{T}.}
    ]{align}
      e^{DA+}_{t,\omega_i}  &=  e^{DA+}_{t,\omega_{i+1}}, 
      &e^{DA-}_{t,\omega_i} &=  e^{DA-}_{t,\omega_{i+1}}    \label{EQ:NAC_e}\\ 
      ie^{DA+}_{t,\omega_i} &=  ie^{DA+}_{t,\omega_{i+1}},
      &ie^{DA-}_{t,\omega_i}&=  ie^{DA-}_{t,\omega_{i+1}}   \label{EQ:NAC_ie}\\ 
      r^{U}_{t,\omega_i}    &=  r^{U}_{t,\omega_{i+1}},
      &r^{D}_{t,\omega_i}   &=  r^{D}_{t,\omega_{i+1}}      \label{EQ:NAC_r}\\ 
      r^{U,B}_{t,\omega_i}  &=  r^{U,B}_{t,\omega_{i+1}}, 
      &r^{D,B}_{t,\omega_i} &=  r^{D,B}_{t,\omega_{i+1}}    \label{EQ:NAC_rB}\\ 
      r^{U,FD}_{t,\omega_i} &=  r^{U,FD}_{t,\omega_{i+1}},
      &r^{D,FD}_{t,\omega_i}&=  r^{D,FD}_{t,\omega_{i+1}}   \label{EQ:NAC_rFD} 
    \end{empheq}
\end{subequations}

For intraday market variables, we have
\begin{equation} \label{EQ:NAC_eIM}
	e^{IM}_{i, t, \omega_{j}} = e^{IM}_{i, t, \omega_{j+1}}, \quad 
    \forall \omega_{j<|\mathcal{C}|_{s,c}} \in \mathcal{C}_{s,c}, \ 
    \forall c \in \mathcal{C}_{S^{IM}_i-1}, \ 
    \forall t \in \mathcal{T}_i, \ 
    \forall i\in \mathcal{I}.
\end{equation}

For flexible demand variables and BESS operational variables, we have
\begin{subequations}
    \begin{empheq}
    [right=\empheqrbrace
    {\ \omega_{i<|\mathcal{C}_{s,c}|}\in\mathcal{C}_{s,c}, c \in \mathcal{C}_{S^{R}_t-1}, t\in\mathcal{T}.}
    ]{align}
      f_{t,\omega_i}       &=  f_{t,\omega_{i+1}} & &   \label{EQ:NAC_f}\\
      {f}^+_{t,\omega_i}  &=  {f}^+_{t,\omega_{i+1}} ,
      &{f}^-_{t,\omega_i} &=  {f}^-_{t,\omega_{i+1}}    \label{EQ:NAC_fpm}\\
      c_{t,\omega_i}        &=  c_{t,\omega_{i+1}},
      &d_{t,\omega_i}       &=  d_{t,\omega_{i+1}},     \label{EQ:NAC_cd}\\
      id_{t,\omega_i}       &=  id_{t,\omega_{i+1}},
      &soc_{t,\omega_i}     &=  soc_{t,\omega_{i+1}}    \label{EQ:NAC_idsoc}
    \end{empheq}
\end{subequations}

Finally, for positive and negative imbalance variables, we have
\begin{equation} \label{EQ:NAC_pPIB}
	e^{IB+}_{t, \omega_i} = e^{IB+}_{t, \omega_{i+1}}, \quad
    e^{IB-}_{t, \omega_i} = e^{IB-}_{t, \omega_{i+1}}, \quad
    \forall \omega_{i<|\mathcal{C}|_{s,c}} \in \mathcal{C}_{s,c}, \ 
    \forall c \in \mathcal{C}_{S_t^R}, \ 
    \forall t \in \mathcal{T}.
\end{equation}

\subsection{Objective Function}
\label{SUBSEC:ObjectiveFunction}

We define the objective function Expected Energy Community Social Welfare (EECSW) as follows:
\begin{subequations}
	\begin{align}
		EECSW & := \sum_{t \in \mathcal{T}} \sum_{\omega \in \Omega} \Pi_{\omega} \lambda^D_{t, \omega} (e^{DA+}_{t, \omega} - e^{DA-}_{t, \omega}) \label{EQ:DAincome} \\
		& + \sum_{t \in \mathcal{T}} \sum_{\omega \in \Omega} \Pi_{\omega} \lambda^R_{t, \omega} (r^D_{t, \omega} + r^U_{t, \omega}) \label{EQ:RESincome} \\
		&  + \sum_{t \in \mathcal{T}} \sum_{\omega \in \Omega} \Pi_{\omega} \sum_{i \in \mathcal{I}_t} \lambda^{I}_{i, t, \omega} e^{IM}_{i, t, \omega} \label{EQ:IMincome} \\
		& + \sum_{t \in \mathcal{T}} \sum_{\omega \in \Omega} \Pi_{\omega} \lambda^{IB+}_{t, \omega} e^{IB+}_{t, \omega} \label{EQ:PosIBcollectRights} \\
		& - \sum_{t \in \mathcal{T}} \sum_{\omega \in \Omega} \Pi_{\omega} \lambda^{IB-}_{t, \omega} e^{IB-}_{t, \omega} \label{EQ:NegIBPayObligations} \\
		& - \sum_{t \in \mathcal{T}} \sum_{\omega \in \Omega} \Pi_{\omega} C^{FD} (f^+_{t, \omega} + f^-_{t, \omega}). \label{EQ:DemFlexPenalty}
	\end{align}
	\label{EQ:EECSW}
\end{subequations}

The first term \eqref{EQ:DAincome} represents the EC's expected income or costs from the day-ahead market. The second term \eqref{EQ:RESincome} represents the EC's expected income from the reserve market. The third term \eqref{EQ:IMincome} represents the EC's expected income or costs from intraday markets. The forth term \eqref{EQ:PosIBcollectRights} represents the EC's expected positive imbalances collection rights. The fifth term \eqref{EQ:NegIBPayObligations} represents the EC's expected imbalance payment obligations. The last term \eqref{EQ:DemFlexPenalty} is a penalty on the use of demand flexibility.

\subsection{Mathematical Formulation Summary} 
\label{SUBSEC:TheCompleteOptimizationModel} 

In this section, we summarise the complete mathematical formulation of the model. Table \ref{TABLE:stages-rv-dv} shows a summary of the stages considered. Each stage can contain state variables, observations and recourse variables. Note that, in terms of the optimization model, recourse variables of one stage are equivalent to decision variables of the next stage.


\begin{table}
\footnotesize
\begin{tabular}{c c c c}
	Stage & State Variables & Observations & Recourse Variables \\
	1 & & $\lambda^D \in \R^{24}$ & $e^{DA+}_{t, \omega}, e^{DA-}_{t, \omega}, ie^{DA+}_{t, \omega}, ie^{DA-}_{t, \omega}$\\
	2 & $r^U_{t, \omega}, r^{U, B}_{t, \omega}, r^{U, FD}_{t, \omega}, r^{D}_{t, \omega}, r^{D, B}_{t, \omega}, r^{D, FD}_{t, \omega}$ & $\lambda^R \in \R^{24}$ & \\
	3 & $e^{IM}_{1, t, \omega}$ & $\lambda^I_1 \in \R^{24}$ &  \\
	4 & $e^{IM}_{2, t, \omega}$ & $\lambda^I_2 \in \R^{24}$ &  \\
	5 & $f_1, f^{\pm}_1, c_{1, \omega}, d_{1, \omega}, id_{1, \omega}, \soc_{1, \omega}$ & $W_1, PV_1, \lambda^{IB\pm}_{1, \omega} \in \R$ & $p^{IB+}_{1, \omega}, p^{IB-}_{2, \omega}$  \\
	6 & $f_2, f^{\pm}_2, c_{2, \omega}, d_{2, \omega}, id_{2, \omega}, \soc_{2, \omega}$ & $W_2, PV_2, \lambda^{IB\pm}_{2, \omega} \in \R$ & $p^{IB+}_{2, \omega}, p^{IB-}_{2, \omega}$ \\
	7 & $e^{IM}_{3, t, \omega}$ & $\lambda^I_3 \in \R^{20}$ &  \\
	8 & $f_3, f^{\pm}_3, c_{3, \omega}, d_{3, \omega}, id_{3, \omega}, \soc_{3, \omega}$  & $W_3, PV_3, \lambda^{IB\pm}_{3, \omega} \in \R$ & $p^{IB+}_{3, \omega}, p^{IB-}_{3, \omega}$ \\
	9 & $f_4, f^{\pm}_4, c_{4, \omega}, d_{4, \omega}, id_{4, \omega}, \soc_{4, \omega}$ & $W_4, PV_4, \lambda^{IB\pm}_{4, \omega} \in \R$ & $p^{IB+}_{4, \omega}, p^{IB-}_{4, \omega}$ \\
	10 & $f_5, f^{\pm}_5, c_{5, \omega}, d_{5, \omega}, id_{5, \omega}, \soc_{5, \omega}$ & $W_5, PV_5, \lambda^{IB\pm}_{5, \omega} \in \R$ & $p^{IB+}_{5, \omega}, p^{IB-}_{5, \omega}$\\
	11 & $e^{IM}_{4, t, \omega}$  & $\lambda^I_4 \in \R^{17}$ &   \\
	12 & $f_6, f^{\pm}_6, c_{6, \omega}, d_{6, \omega}, id_{6, \omega}, \soc_{6, \omega}$ & $W_6, PV_6, \lambda^{IB\pm}_{6, \omega} \in \R$ & $p^{IB+}_{6, \omega}, p^{IB-}_{6, \omega}$ \\
	13 & $f_7, f^{\pm}_7, c_{7, \omega}, d_{7, \omega}, id_{7, \omega}, \soc_{7, \omega}$ & $W_7, PV_7, \lambda^{IB+}_{7, \omega} \in \R$ & $p^{IB+}_{7, \omega}, p^{IB-}_{7, \omega}$ \\
	14 & $f_8, f^{\pm}_8, c_{8, \omega}, d_{8, \omega}, id_{8, \omega}, \soc_{8, \omega}$ & $W_8, PV_8, \lambda^{IB+}_{8, \omega} \in \R$ & $p^{IB+}_{8, \omega}, p^{IB-}_{8, \omega}$ \\
	15 & $f_9, f^{\pm}_9, c_{9, \omega}, d_{9, \omega}, id_{9, \omega}, \soc_{9, \omega}$ & $W_9, PV_9, \lambda^{IB+}_{9, \omega} \in \R$ & $p^{IB+}_{9, \omega}, p^{IB-}_{9, \omega}$\\
	16 & $e^{IM}_{5, t, \omega}$ & $\lambda^I_5 \in \R^{13}$ &   \\
	17 & $f_{10}, f^{\pm}_{10}, c_{10, \omega}, d_{10, \omega}, id_{10, \omega}, \soc_{10, \omega}$ & $W_{10}, PV_{10}, \lambda^{IB\pm}_{10, \omega} \in \R$ & $p^{IB+}_{10, \omega}, p^{IB-}_{10, \omega}$ \\
	18 & $f_{11}, f^{\pm}_{11}, c_{11, \omega}, d_{11, \omega}, id_{11, \omega}, \soc_{11, \omega}$ & $W_{11}, PV_{11}, \lambda^{IB\pm}_{11, \omega} \in \R$ & $p^{IB+}_{11, \omega}, p^{IB-}_{11, \omega}$ \\
	19 & $f_{12}, f^{\pm}_{12}, c_{12, \omega}, d_{12, \omega}, id_{12, \omega}, \soc_{12, \omega}$ & $W_{12}, PV_{12}, \lambda^{IB\pm}_{12, \omega} \in \R$ & $p^{IB+}_{12, \omega}, p^{IB-}_{12, \omega}$ \\
	20 & $f_{13}, f^{\pm}_{13}, c_{13, \omega}, d_{13, \omega}, id_{13, \omega}, \soc_{13, \omega}$ & $W_{13}, PV_{13}, \lambda^{IB\pm}_{13, \omega} \in \R$ & $p^{IB+}_{13, \omega}, p^{IB-}_{13, \omega}$ \\
	21 & $e^{IM}_{6, t, \omega}$ & $\lambda^I_6 \in \R^{9}$ &   \\
	22 & $f_{14}, f^{\pm}_{14}, c_{14, \omega}, d_{14, \omega}, id_{14, \omega}, \soc_{14, \omega}$ & $W_{14}, PV_{14}, \lambda^{IB\pm}_{14, \omega} \in \R$ & $p^{IB+}_{14, \omega}, p^{IB-}_{14, \omega}$ \\
	23 & $f_{15}, f^{\pm}_{15}, c_{15, \omega}, d_{15, \omega}, id_{15, \omega}, \soc_{15, \omega}$ & $W_{15}, PV_{15}, \lambda^{IB\pm}_{15, \omega} \in \R$ & $p^{IB+}_{15, \omega}, p^{IB-}_{15, \omega}$ \\
	24 & $f_{16}, f^{\pm}_{16}, c_{16, \omega}, d_{16, \omega}, id_{16, \omega}, \soc_{16, \omega}$ & $W_{16}, PV_{16}, \lambda^{IB\pm}_{16, \omega} \in \R$ & $p^{IB+}_{16, \omega}, p^{IB-}_{16, \omega}$ \\
	25 & $f_{17}, f^{\pm}_{17}, c_{17, \omega}, d_{17, \omega}, id_{17, \omega}, \soc_{17, \omega}$ & $W_{17}, PV_{17}, \lambda^{IB\pm}_{17, \omega} \in \R$ & $p^{IB+}_{17, \omega}, p^{IB-}_{17, \omega}$ \\
	26 & $f_{18}, f^{\pm}_{18}, c_{18, \omega}, d_{18, \omega}, id_{18, \omega}, \soc_{18, \omega}$ & $W_{18}, PV_{18}, \lambda^{IB\pm}_{18, \omega} \in \R$ &  $p^{IB+}_{18, \omega}, p^{IB-}_{18, \omega}$ \\
	27 & $f_{19}, f^{\pm}_{19}, c_{19, \omega}, d_{19, \omega}, id_{19, \omega}, \soc_{19, \omega}$ & $W_{19}, PV_{19}, \lambda^{IB\pm}_{19, \omega} \in \R$ & $p^{IB+}_{19, \omega}, p^{IB-}_{19, \omega}$ \\
	28 & $e^{IM}_{7, t, \omega}$ & $\lambda^I_7 \in \R^{4}$ &   \\
	29 & $f_{20}, f^{\pm}_{20}, c_{20, \omega}, d_{20, \omega}, id_{20, \omega}, \soc_{20, \omega}$  & $W_{20}, PV_{20}, \lambda^{IB\pm}_{20, \omega} \in \R$ & $p^{IB+}_{20, \omega}, p^{IB-}_{20, \omega}$ \\
	30 &  $f_{21}, f^{\pm}_{21}, c_{21, \omega}, d_{21, \omega}, id_{21, \omega}, \soc_{21, \omega}$ & $W_{21}, PV_{21}, \lambda^{IB\pm}_{21, \omega} \in \R$ & $p^{IB+}_{21, \omega}, p^{IB-}_{21, \omega}$ \\
	31 & $f_{22}, f^{\pm}_{22}, c_{22, \omega}, d_{22, \omega}, id_{22, \omega}, \soc_{22, \omega}$ & $W_{22}, PV_{22}, \lambda^{IB\pm}_{22, \omega} \in \R$ & $p^{IB+}_{22, \omega}, p^{IB-}_{22, \omega}$ \\
	32 & $f_{23}, f^{\pm}_{23}, c_{23, \omega}, d_{23, \omega}, id_{23, \omega}, \soc_{23, \omega}$ & $W_{23}, PV_{23}, \lambda^{IB\pm}_{23, \omega} \in \R$ & $p^{IB+}_{23, \omega}, p^{IB-}_{23, \omega}$ \\
	33 & $f_{24}, f^{\pm}_{24}, c_{24, \omega}, d_{24, \omega}, id_{24, \omega}, \soc_{24, \omega}$ & $W_{24}, PV_{24}, \lambda^{IB\pm}_{24, \omega} \in \R$ & $p^{IB+}_{24, \omega}, p^{IB-}_{24, \omega}$ \\
\end{tabular}
\caption{Summary of stages, decision variables, random variable observations and recourse variables.}
\label{TABLE:stages-rv-dv}
\end{table}

Finally, we present a summary of the complete mathematical formulation.

\begin{align*}
\notag & \text{max} & & EECSW \\
\notag & \text{s.t.} & & \eqref{EQ:FlexDemandPos}, \eqref{EQ:DailyDem}, \eqref{EQ:InterDem}, \eqref{EQ:FlexDemDisplace},  & \text{Flex. Demand Op.} \\
\notag & & & \eqref{dVdCPos}-\eqref{EQ:SOCV_fin}, & \text{BESS Op.} \\
\notag & & & \eqref{EQ:DAVarDef}-\eqref{EQ:DA_bid_mono_2}, & \text{EC DA} \\
\notag & & & \eqref{EQ:RVarDef}-\eqref{EQ:VPP_RD_SOCV},  & \text{EC RM} \\
\notag & & & \eqref{EQ:DM_IM_1}-\eqref{EQ:IM_bounds_3-4}, & \text{EC IM} \\
\notag & & & \eqref{EQ:IBVarDef}-\eqref{EQ:VPP_IB_MIN}, & \text{EC IB} \\
\notag & & & \eqref{EQ:NAC_e}-\eqref{EQ:NAC_idsoc}. & \text{NAC} \\
\end{align*}


\bibliographystyle{elsarticle-num} 
\bibliography{references}






\end{document}